\documentclass[reqno,11pt,noamsfonts]{amsart}
\usepackage[foot]{amsaddr}         

\usepackage[default,scale=0.92]{opensans}
\usepackage[notext]{stix}

\usepackage{etoolbox}
\patchcmd{\section}{\scshape}{\bfseries}{}{}
\makeatletter
\renewcommand{\@secnumfont}{\bfseries}
\makeatother

\usepackage{amsthm}
\newtheorem{remark}{Remark}
\newtheorem{proposition}{Proposition}

\usepackage{bm}
\usepackage{mathtools}
\usepackage{xfrac}                  % for \sfrac
\usepackage{accents}                % for \underaccent

\usepackage[svgnames]{xcolor}

\usepackage[numbers,sort&compress]{natbib} 

\usepackage[numbers,sort&compress]{natbib} 

\usepackage{graphicx}
\usepackage[labelfont=bf]{caption}
\usepackage{subcaption}
\usepackage{enumitem}
\usepackage{array}             % use these 
\usepackage{booktabs}          % for tables

\definecolor{Mosss}{RGB}{0,192,81}
\definecolor{Reddd}{RGB}{240,40,40}
\definecolor{Mocca}{RGB}{148,82,0}

\calclayout

\renewcommand {\mathit}[1]{\text{{\rmfamily\itshape #1}}}

\newcommand {\CFL}    {\mathit{CFL}}

\newcommand {\EX}     {\mathrm{ex}}       % "explicit" part
\newcommand {\IM}     {\mathrm{im}}       % "implicit" part

\newcommand {\fc}     {\M f_\mathrm{c}}

\newcommand {\fd}     {\M f_\mathrm{d}}
\newcommand {\fs}     {\M f_\mathrm{s}}
\newcommand {\fex}    {\M \phi_\mathrm{ex}}
\newcommand {\fim}    {\M \phi_\mathrm{im}}

\newcommand {\Ac}     {\M A_\mathrm{c}}
\newcommand {\Ad}     {\M A_{\mathrm{d}}}

\newcommand {\PO}  {\mathcal M}
\newcommand {\F}   {\mathcal F}
\newcommand {\Fc}  {\mathcal F_{\!\mathrm{c}}}
\newcommand {\Fd}  {\mathcal F_{\!\mathrm{d}}}
\newcommand {\Fs}  {\mathcal F_{\!\mathrm{s}}}
\newcommand {\Fex} {\M \varPhi_{\mathrm{ex}}}
\newcommand {\Fim} {\M \varPhi_{\mathrm{im}}}

\newcommand {\iu}  {\mathrm{i}}

\newcommand {\M}  [1] {\bm{#1}}                        % matrix
\newcommand {\MC} [1] {\mathbf{#1}}                    % matrix constant

\newcommand {\NM} [1] {\underaccent{\bar}{#1}}         % coefficient matrix (general or time)
\providecommand {\norm}[1] {\left\Vert#1\right\Vert}

\renewcommand {\d}            {\partial}
\newcommand   {\D}            {\mathrm{d}\mspace{1.0mu}}

\newcommand   {\transpose}[1] {#1^{\mathrm{t}}}
\newcommand   {\jmp}      [1] {\left\lBrack #1\right\rBrack}   % with stix
\newcommand   {\avg}      [1] {\left\lBrace #1\right\rBrace}   % with stix 

\newcommand   {\et}   {E}               % total energy per unit mass.
\newcommand   {\htot} {H}               % total enthalpy per unit mass

\begin{document}

\title%[]
{% Full title
Stable semi-implicit SDC methods for conservation laws
}
\author{Jörg Stiller}
\address{%
  TU Dresden, Institute of Fluid Mechanics, 01062 Dresden, Germany}
\email{joerg.stiller@tu-dresden.de}
%\email{\{joerg.stiller,\dots\}@tu-dresden.de}

\begin{abstract}
%%%
Semi-implicit spectral deferred correction (SDC) methods provide a systematic approach  to construct time integration methods of arbitrarily high order for nonlinear evolution equations including conservation laws.
%%%
They converge towards $A$- or even $L$-stable collocation methods, but are often not sufficiently robust themselves.
%%%
In this paper, a family of SDC methods inspired by an implicit formulation of the Lax-Wendroff method is developed.
%%%
Compared to fully implicit approaches, the methods have the advantage that they only require the solution of positive definite or semi-definite linear systems.
%%%
Numerical evidence suggests that the proposed semi-implicit SDC methods with Radau points are $L$-stable up to order 11 and require very little diffusion for orders 13 and 15.
%%%
The excellent stability and accuracy of these methods is confirmed by numerical experiments with 1D conservation problems, including the convection-diffusion, Burgers, Euler and Navier-Stokes equations.
%%%
\end{abstract}

\keywords{%
Semi-implicit methods; 
Spectral deferred correction;
discontinuous Galerkin method
}

\maketitle

%%%%%%%%%%%%%%%%%%%%%%%%%%%%%%%%%%%%%%%%%%%%%%%%%%%%%%%%%%%%%%%%%%%%%%%%%%%%%%%%%%%%
% Text

% !TEX encoding = UTF-8 Unicode
% !TEX root =  ssi-sdc-1d.tex

%%%%%%%%%%%%%%%%%%%%%%%%%%%%%%%%%%%%%%%%%%%%%%%%%%%%%%%%%%%%%%%%%%%%%%%%%%%%%%%%%%%%

\section{Introduction}
\label{sec:introduction}

%%%
The trend towards higher-order numerical methods for conservation laws in fluid mechanics and beyond increases the need for time integration methods that match the convergence and efficiency of spatial methods.
%%%
A wide range of those methods have been developed in recent decades \cite{TI_Hairer1993a,TI_Hairer1996a,TI_Gottlieb2016a}.
%%%
They include explicit, implicit and semi-implicit methods using multiple steps, stages or derivatives \cite{SE_Canuto2007a}.
%%%
Convection-dominant problems are often not stiff and can be solved with explicit methods.
%%%
In particular, strong stability preserving Runge-Kutta methods are a popular choice for hyperbolic problems with discontinuous solutions \cite{TI_Gottlieb2005a,TI_Ketcheson2008a}.
%%%
However, explicit methods lose efficiency when the stiffness increases due to unresolved waves, diffusion or local mesh refinement.
%%%
This issue is particularly pronounced with high-order spatial finite element methods, where the diffusion time scales like ${P^4/h^2}$ for elements of size $h$ and degree $P$ \cite{SE_Canuto2011a}.
%%%
One approach to avoid the resulting time step restrictions is to use fully implicit methods such as 
diagonally implicit Runge-Kutta methods \cite{TI_Kennedy2016a,TI_Pazner2017a,TI_Pan2021a}, 
Rosenbrock methods 
\cite{TI_Bassi2015a},
Taylor methods
\cite{TI_Baeza2020a}
or discontinuous Galerkin methods in time 
\cite{TI_Tavelli2018a}.
%%%
Unfortunately, these methods often lead to ill conditioned algebraic systems that are expensive to solve.
%%%
In time-resolved simulations, it therefore remains a challenge to compete with explicit methods.
%%%

Semi-implicit methods take advantage of the fact that the asymmetric and often nonlinear convection terms are non-stiff and thus suited for explicit treatment. 
%%%
The diffusion terms usually require an implicit discretization, but are (quasi-)linear and lead to (semi-)definite algebraic systems that can be solved with standard methods of linear algebra.
%%%
Popular semi-implicit methods include implicit-explicit (IMEX) multistep methods and Runge-Kutta methods.
%%%
IMEX multistep methods can  achieve an arbitrary order of convergence by extrapolating the non-stiff terms \cite[Ch.4.2]{TI_Hundsdorfer2003a}.
%%%
However, the increase in order affects stability \cite{TI_Ascher1995a,TI_Frank1997a}.
%%%
Therefore, the IMEX multistep methods used in practice rarely exceed order 3, e.g.
 \cite{TI_Karniadakis1991a,TI_Fehn2017a,TI_Klein2015a}.
%%%
IMEX Runge-Kutta methods are often a better choice because they possess a larger stability domain and a lower error constant than comparable multistep methods \cite{TI_Ascher1997a,TI_Kennedy2003a,TI_Cavaglieri2015a,TI_Izzo2017a}.
%%%
By using them as a partitioned scheme, IMEX Runge methods also allow for the
semi-implicit treatment of quasi-linear stiff terms \cite{TI_Boscarino2016a}.
This approach was recently applied to compressible and incompressible Navier-Stokes problems \cite{TI_Boscheri2022a,TI_Guesmi2023a}.
%%%
However, due to the increasingly complex order conditions, it is difficult to increase the rate of convergence and there is no known IMEX Runge-Kutta method with an order higher than 5.

%%%
Recently, considerable progress in developing higher order semi-implicit time integration methods was achieved using 
multi-derivative methods, e.g.
\cite{TI_Schuetz2021a,TI_Frolkovic2023a,TI_Macca2024a}.
%%%
Alternatively, spectral deferred correction (SDC) methods can be used to construct schemes of arbitrary high order \cite{TI_Dutt2000a}.
%%%
Semi-implicit SDC methods were proposed by \citet{TI_Minion2003b} and further developed in \cite{TI_Hagstrom2006a,TI_Layton2007a,TI_Christlieb2011a,TI_Christlieb2015a}.
%%%
In recent applications to incompressible Navier-Stokes problems, these methods achieved temporal convergence rates up to 12 \cite{TI_Minion2018a,TI_Stiller2020a}.
%%%
For laminar flow problems, semi-implicit SDC methods proved competitive with IMEX multistep and Runge-Kutta methods for relative $L^2$ errors below $10^{-4}$.
%%%
However, in direct simulations of turbulent flows, where non-linear convective processes dominate, the SDC methods suffered from instabilities \cite{TI_Guesmi2023a}.
%%%

For improving the usability of semi-implicit SDC methods it is imperative to increase their stability.
%%%
To achieve this goal, a novel approach is proposed in this paper.
%%%
It divides the convection term into a possibly nonlinear explicit part and an implicitly treated diffusion-like part.
%%%
This splitting is inspired by the Lax-Wendroff method, where the two parts appear as the first terms of a Taylor (or Cauchy-Kovalevskaya) expansion \cite{TI_Lax1960a}.
%%%
The main difference consists in the implicit treatment of the diffusion-like part.
%%%
A similar strategy was already used in Taylor-Galerkin methods \cite{TI_Donea1984a,TI_Safjan1995a}, streamline-upwind Petrov-Galerkin methods \cite{FE_Brooks1982a,FE_Hughes1986c} and multi-derivative methods \cite{TI_Jaust2016a,TI_Frolkovic2023a}.
%%%
In contrast to these methods, however, the present work does not extend the Taylor series expansion to the physical diffusion terms. 
%%%
Instead, the incomplete expansion is used to design $L$-stable time integrators, from which robust high-order SDC methods are constructed.
%%%
Due to the different treatment of the convection term, these SDC methods are three-quarter-implicit rather than semi-implicit.
%%%
Nevertheless, they require only the solution of well-posed linear algebraic equations and thus share the advantages of conventional semi-implicit methods.

%%%
An outline of the paper follows: 
%%%
Section~\ref{sec:conservation-laws} gives a brief overview of the one-dimensional conservation laws considered.
%%%
The first part of Section~\ref{sec:time-integration} is concerned with the development of 
robust semi-implicit one- and two-stage integrators. In the second part, SDC methods based on these integrators are constructed and analyzed.
%%%
Spatial discretization including stabilization techniques is described in Section~\ref{sec:spatial-discretization}, followed by a brief discussion of solution methods in Section~\ref{sec:solution}.
%%%
Section~\ref{sec:numerical-experiments} provides numerical experiments for convection-diffusion, Burgers, and compressible flow problems.
%%%
Conclusions are drawn in Section~\ref{sec:conclusions}.
%%%

% !TEX encoding = UTF-8 Unicode
% !TEX root =  ssi-sdc-1d.tex

%%%%%%%%%%%%%%%%%%%%%%%%%%%%%%%%%%%%%%%%%%%%%%%%%%%%%%%%%%%%%%%%%%%%%%%%%%%%%%%%%%%%

\section{Conservation laws}
\label{sec:conservation-laws}

%===================================================================================

\subsection{Conservation systems in one space dimension}
\label{sec:conservation-laws:systems}

%-----------------------------------------------------------------------------------

Consider the one-dimensional conservation law
\begin{equation}
  \label{eq:conservation-system}
  \d_t \M u = -\d_x \fc(\M u) + \d_x \fd(\M u) + \fs(\M u, x,t)
            \eqqcolon \M f(\M u, x,t)
\end{equation}
for 
  ${x \in \Omega \subset \mathbb R}$
and
  ${t \in \mathbb R^{+}}$,
where
  ${\M u(x,t) \in \mathbb R^d}$ 
are the conservative variables,
  $\fc$
the convective fluxes, 
  $\fd$
the diffusive fluxes,
  $\fs$
the sources and 
  $\M f$
is the combined right-hand side.
%%%
It is assumed that the convective fluxes possess the Jacobian ${\Ac = \fc'(\M u)}$ and
the diffusive fluxes can be written in the form ${\fd = \Ad \d_x \M u}$, where
 ${\Ad(\M u,\d_x \M u)}$
is the positive-semidefinite diffusion matrix.
%%%
These definitions yield the quasilinear form
\begin{equation*}
  %\label{eq:conservation-system:quasilinear}
  \d_t \M u = -\Ac \d_x \M u + \d_x \Ad \d_x \M u + \fs
  \,.
\end{equation*}

%===================================================================================

\subsection{Scalar conservation laws}
\label{sec:conservation-laws:scalar}

%-----------------------------------------------------------------------------------

Scalar laws consist of only one conservation equation with solution $u$ such that $\M u = [u]$.
%%%
Two cases are considered.
%%%
The first is the convection-diffusion equation
\begin{equation}
  \label{eq:conv-diff}
  \d_t u = -\d_x (vu) + \d_x(\nu \d_x u)
\end{equation}
with constant velocity $v$ and diffusivity ${\nu \ge 0}$.
%%%
This yields the convective Jacobian
  ${\Ac = [v]}$
and the diffusion matrix
  ${\Ad = [\nu]}$.
%%%
The second case is the Burgers equation
\begin{equation*}
  %\label{eq:burgers}
  \d_t u = - \d_x \bigl(\tfrac{1}{2}u^2\bigr) 
           + \d_x(\nu \d_x u) +  f_{\mathrm s}(x,t)
\end{equation*}
with constant diffusivity ${\nu \ge 0}$,
convective Jacobian
  ${\Ac = [u]}$
and diffusion matrix
  ${\Ad = [\nu]}$.
%%%
The source $f_{\mathrm s}$ can be manufactured to match a prescribed exact solution.
%%%

%===================================================================================

\subsection{Euler and  Navier-Stokes equations}
\label{sec:conservation-laws:cns}

For the Navier-Stokes equations, the conservative variables, convective fluxes and diffusive fluxes given by
\begin{equation*}
  \M u = \begin{bmatrix}
           \rho \\ \rho v \\ \rho \et
         \end{bmatrix}
  \,, \quad
  \fc = \begin{bmatrix}
          \rho v \\ \rho v^2 + p\\ \rho v \et + pv
        \end{bmatrix}
  \,, \quad
  \fd = \begin{bmatrix}
           0 \\ 
           \frac{4}{3} \eta \d_x v \\ 
           \frac{4}{3} \eta v \d_x v + \lambda \d_x T
        \end{bmatrix}
  \,, 
%  , \quad
%  \fs = \begin{bmatrix}
%           0 \\ f(x,t) \\ v f(x,t)
%        \end{bmatrix}    
\end{equation*}
where
  $\rho$ 
is the density,
  $v$
the velocity,
  $\et$ 
the total specific energy,
  $\eta$
the dynamic viscosity and
  $\lambda$
the heat conductivity.
%%%
Assuming a perfect gas with the gas constant ${R = c_p - c_v}$ and the ratio of the specific heats ${\gamma = c_p/c_v}$ it is possible to calculate the temperature
  ${T = (\et - v^2/2)/c_v}$,
the pressure
  ${p = \rho R T}$,
the total specific enthalpy ${\htot = \et + p/\rho}$
and
the speed of sound ${a = ((\gamma-1) R T)^{1/2}}$.
%%%
The convective Jacobian and the diffusion matrix are given by
\begin{equation*}
  %\label{eq:convective-jacobian:cns}
  \Ac =
      \begin{bmatrix}
        0                                 & 1                     & 0        \\[1mm]
        \frac{1}{2}(\gamma-3)v^2          & (3-\gamma)v           & \gamma-1 \\[1mm]
        \frac{1}{2}(\gamma-1)v^3 - v\htot & \htot - (\gamma-1)v^2 & \gamma v
      \end{bmatrix}
\end{equation*}
and
\begin{equation*}
  %\label{eq:diffusion-matrix:cns}
  \Ad =
     \begin{bmatrix}
          0                                               
       &  0                      
       &  0 
       \\[1mm]
         -\frac{4}{3}\nu v                              
       &  \frac{4}{3}\nu         
       &  0
       \\[1mm]
         -\bigl(\frac{4}{3}\nu - \gamma\kappa\bigr)v^2 - \gamma\kappa \et
       &  \bigl(\frac{4}{3}\nu - \gamma\kappa\bigr)v 
       &  \gamma\kappa
     \end{bmatrix}
     \,,
  \end{equation*}
where
${\nu = \eta/\rho}$
is the kinematic viscosity
and
${\kappa = \lambda /(\rho c_p)}$
the thermal diffusivity.
%%%
The inviscid case with ${\eta = \lambda = 0}$ leads to the Euler equations for which
${\Ad = \M 0}$.

%===================================================================================

% !TEX encoding = UTF-8 Unicode
% !TEX root = ssi-sdc-1d.tex

%%%%%%%%%%%%%%%%%%%%%%%%%%%%%%%%%%%%%%%%%%%%%%%%%%%%%%%%%%%%%%%%%%%%%%%%%%%%%%%%%%%%

\section{Time integration}
\label{sec:time-integration}

%===================================================================================

\subsection{Semi-explicit integrators of order 1 and 2}
\label{sec:time-integration:integrators}

This section provides low-order time integrators as building blocks for SDC methods.
%%%
The aim is to devise a semi-implicit approach that avoids the fully implicit treatment of convection, but achieves A-stability and exhibits a robust behavior when applied to nonlinear conservation laws.
%%%
For readability, the dependence on $x$ is hidden throughout the section, except for the spatial derivatives.

As the starting point consider the purely convective system
\begin{equation}
  \label{eq:conservation-system:convective}
  \d_t \M u + \d_x \fc(\M u) = 0
  \,.
\end{equation}
The application of the Lax-Wendroff method \cite{TI_Lax1960a} to \eqref{eq:conservation-system:convective} yields
\begin{equation*}
  %\label{eq:lax-wendroff}
  \frac{\M u^{n+1} - \M u^{n}}{\Delta t} 
  + \d_x \fc(\M u^{n}) 
  = \d_x \frac{\Delta t}{2} \Ac^2(\M u^{n})\,\d_x \M u^{n}
  \,,
\end{equation*}
where
${\Delta t = t^{n+1} - t^{n}}$ is the time step size and
${\M u^n \approx \M u(t^{n})}$.
%%%
This scheme is second-order accurate but only conditionally stable.
%%%
To improve the stability, the right-hand side can be changed to a semi-implicit form such that
\begin{equation}
  \label{eq:lax-wendroff:modified}
  \frac{\M u^{n+1} - \M u^{n}}{\Delta t} 
  + \d_x \fc(\M u^{n})
  = \d_x \frac{\Delta t}{2} \Ac^2(\M u^{n})\,\d_x \M u^{n+1}
  \,.
\end{equation}
The modified scheme retains order 2 and proves unconditionally stable for linear convection with constant velocity (see Sec.~\ref{sec:time-integration:integrators:stability}).
%%%
Its implicit part corresponds to a diffusion problem, which yields a positive-definite  linear system when discretized in space.
%%%
Therefore, the discrete equations are significantly easier to solve than the nonlinear or indefinite systems that emerge from fully implicit methods.
%%%

Unfortunately, there exists no straightforward extension of the Lax-Wendroff method to conservation laws including diffusion.
%%%
Constructing a corresponding single stage method of order 2 would imply the introduction of third derivatives as in Taylor-Galerkin methods or multiderivative methods, see e.g. \cite{TI_Donea1984a,TI_Hairer1993a}.
%%%
The methods proposed below avoid this complication by introducing further stages.

%-----------------------------------------------------------------------------------

\subsubsection{Integrators of order 1}
\label{sec:time-integration:integrators:order-1}

For concise notation, the functionals
\begin{subequations}
  \label{eq:rhs:imex:semi}
  \begin{align}
    \fex\big(\M u^\alpha\big)
    & = -\d_x \fc\big(\M u^\alpha\big)
    \,,
    \\
    \fim\big(\M u^\alpha, \M u^\beta, t, \theta\big)
    & = \d_x \Bigl( \Bigl( \frac{\theta}{2} \Ac^2\big(\M u^\alpha\big) 
                         + \Ad\big(\M u^\alpha\big) 
                     \Bigr) \d_x\M u^\beta
              \Bigr)
      + \fs\big(\M u^\alpha,x,t\big)
  \end{align}
\end{subequations}
are introduced to distinguish the explicit and implicit contributions of convection, diffusion and sources to the semi-discrete equations. 
In \eqref{eq:rhs:imex:semi} the arguments ${\M u^\alpha(x)}$ and ${\M u^\beta(x)}$ represent different instances of the semi-discrete solution or approximations to it, while $\theta$ is a parameter which corresponds to an interval in time. 
The explicit dependence of $\fim$ on $x$ is dropped for brevity.

Using these definitions and accounting for ${\Ad = \MC 0}$, the modified Lax-Wendroff method \eqref{eq:lax-wendroff:modified} can be rewritten in the form
\begin{equation}
  \label{eq:SI1(1):semi}
  \M u^{n+1} 
  = \M u^{n}
  + \Delta t \bigl[ \fex(\M u^{n})
                  + \fim(\M u^n, \M u^{n+1}, t^{n+1}, \Delta t)
             \bigr]
  \,.
\end{equation}
Given the above definition for $\fim$, this scheme applies also to convective-diffusive conservation laws \eqref{eq:conservation-system}, for which it yields a backward Euler discretization of the diffusion term.
Sources are treated implicitly if they are independent of the solution and explicitly otherwise.
In summary, \eqref{eq:SI1(1):semi} can be characterized one-stage semi-implicit time integrator of order 1.
For simplicity, the method is referred to as SI1(1) in the following. 

%%%
Although SI1(1) has excellent stability properties (see Sec.~\ref{sec:time-integration:integrators:stability}), a stronger damping was found to be necessary in high-order SDC methods.
%%%
This is achieved by a second stage in which the result of \eqref{eq:SI1(1):semi} is used to update the explicit part of the RHS.
%%%
The resulting two-stage first-order method SI1(2) is given by
\begin{equation}
  \label{eq:SI1(2):semi}
  \begin{alignedat}{2}
    & \M u^{(1)} 
    && = \M u^{n}
       + \Delta t \bigl[ \fex(\M u^{n})
                       + \fim(\M u^n, \M u^{(1)}, t^{n+1}, \Delta t)
                  \bigr]
    \,,
    \\
    & \M u^{(2)} 
    && = \M u^{n}
       + \Delta t \bigl[ \fex(\M u^{(1)})
                       + \fim( \M u^n, \M u^{(2)},t^{n+1}, \Delta t)
                  \bigr]
    \,,
    \\
    &\M u^{n+1} &&= \M u^{(2)}
    \,.
  \end{alignedat}
\end{equation}
%%%
Using Taylor expansion, it can be shown that the second stage adds a term of the form ${\d_x (\Delta t \Ac^2 \d_x\M u)}$, which is about twice the amount of numerical diffusion introduced by the Euler backward method.
%%%
A similar effect could be achieved without a second stage by including an artificial diffusion term in SI1(1).
%%%
However, comparative studies have shown that this approach is less efficient.

%-----------------------------------------------------------------------------------

\subsubsection{Integrator of order 2}
\label{sec:time-integration:integrator:order-2}

To obtain a method of order 2, SI1(2) can be modified by performing two half steps followed by the application of the midpoint rule.
%%%
This yields the semi-implicit two-stage method SI2(2)
\begin{equation*}
  %\label{eq:SI2(2):semi}
  \begin{alignedat}{2}
    & \M u^{(1)} 
    && = \M u^{n}
       + \frac{\Delta t}{2} 
         \bigl[ \fex(\M u^{n})
              + \fim(\M u^n, \M u^{(1)}, t^n + \Delta t/2, \Delta t)
         \bigr]
    \,,
    \\
    & \M u^{(2)} 
    && = \M u^{n}
       + \frac{\Delta t}{2} 
         \bigl[ \fex(\M u^{(1)})
              + \fim(\M u^n, \M u^{(2)}, t^n + \Delta t/2, \Delta t)
                  \bigr]
    \,,
    \\
    & \M u^{n+1} 
    && = \M u^{n} 
       + \Delta t \M f(\M u^{(2)}, t^n + \Delta t/2)
    \,.
  \end{alignedat}
\end{equation*}
%%%
Another two-stage second-order method can be constructed using the trapezoidal rule.
%%%
However, this method is less stable and therefore not pursued any further.

%-----------------------------------------------------------------------------------

\subsubsection{Linear stability of semi-implicit integrators}
\label{sec:time-integration:integrators:stability}

Consider the convection-diffusion equation \eqref{eq:conv-diff} with the spatially periodic complex solution ${u = \hat u(t) \mathrm{e}^{\iu k x}}$.
%%%
Evaluating the spatial derivatives and removing the common factor $\mathrm{e}^{\iu k x}$ yields
\begin{equation}
   \label{eq:conv-diff:model}
   \D_t \hat u = -\iu v k \hat u  - \nu k^2 \hat u
  \,.
\end{equation}
%%%
This equation is equivalent to the Dahlquist equation
\begin{equation}
  \label{eq:dahlquist}
  \D_t \hat u = \lambda \hat u
  \,,
\end{equation}
where 
${\lambda \in \mathbb C}$
with the real part
${\lambda_{\mathrm r} = -\nu k^2}$
and imaginary part
${\lambda_{\mathrm i} = -v k}$.
%%%
The application of a one-step method to \eqref{eq:dahlquist} yields the stability function
\begin{equation*}
  R(z) = \frac{\hat u^{n+1}}{\hat u^n}
  \,,
\end{equation*}
where ${z = \Delta t \lambda}$.
%%%
The set ${\{z \in \mathbb C; |R(z)| \le 1\}}$ is the stability domain of the method. 
%%%
The method is called 
$A$-stable if ${|R(z)| \le 1}$ for all $z$ with ${z_{\mathrm i} \le 0}$ and
$L$-stable if in addition ${R = 0}$ for ${z \rightarrow \infty}$.
%%%

\begin{remark}
For conservation laws, 
${|z_i| = \Delta t |v| k}$ can be associated with the CFL number and
${-z_r = \Delta t \nu k^2}$ with the diffusion number 
of the fully discrete problem.
\end{remark}

\begin{figure}
%  \subcaptionbox{SI1(1)
%    \label{fig:stability-function:si1}}
    {\includegraphics[scale=0.52]{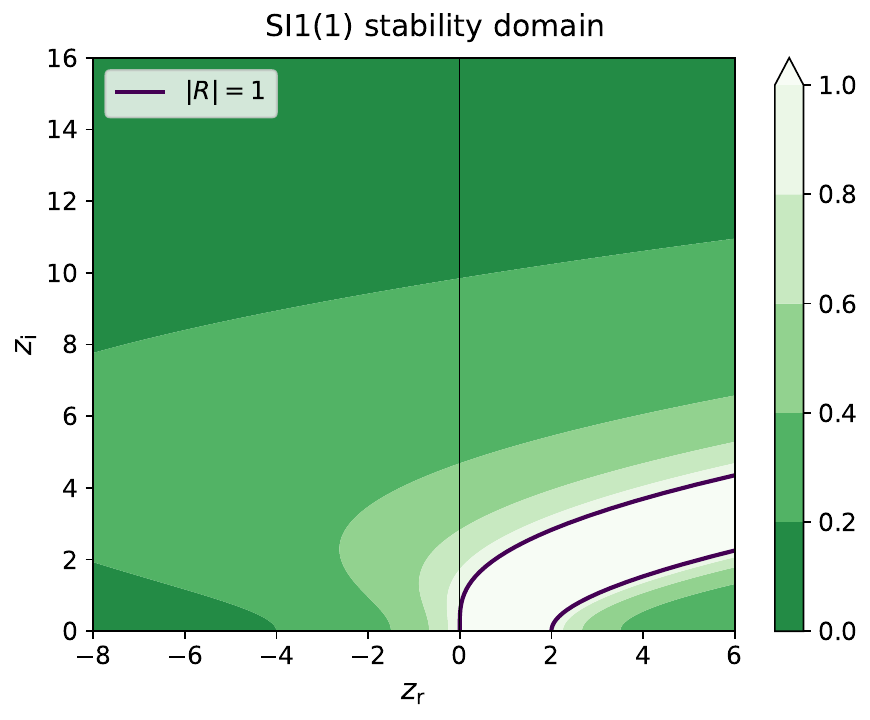}}
%  \subcaptionbox{SI1(2)
%    \label{fig:stability-function:sdc:si12}}
    {\includegraphics[scale=0.52]{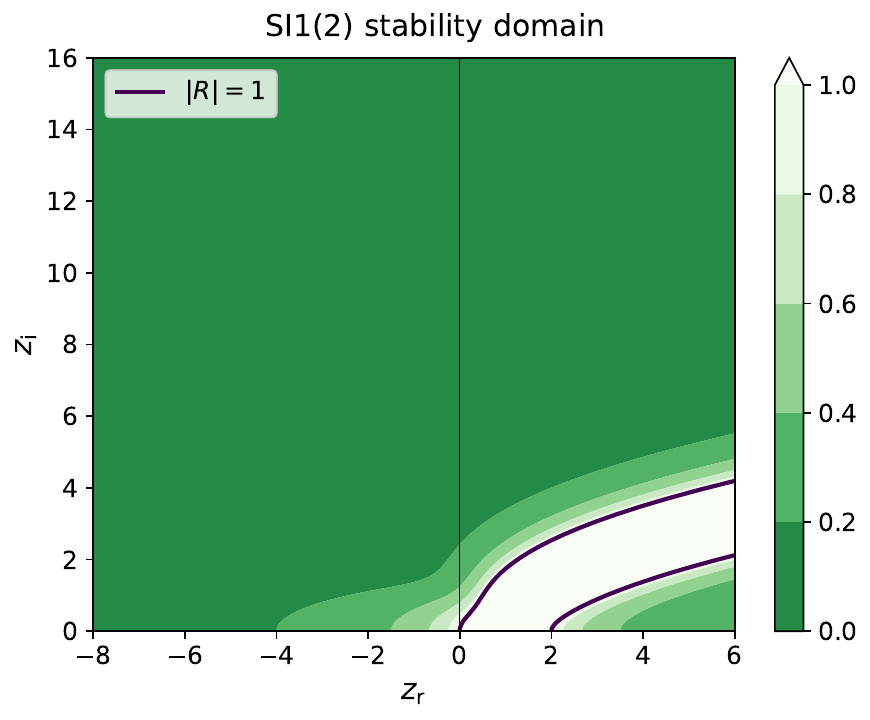}} 
  \caption{Stability domains of the first-order integrators
    \label{fig:stability-functions:si1}}
\end{figure}

To apply the one-step methods presented above to the model problem \eqref{eq:conv-diff:model}, it is necessary to identify the amplitudes of the explicit and implicit RHS parts defined in \eqref{eq:rhs:imex:semi}.
Substituting the solution, evaluating the derivative and separating the factor $\mathrm{e}^{\iu k x}$ gives
%
%\begin{subequations}
  %\label{eq:rhs:imex:conv-diff:model}
  \begin{align*}
    \hat \phi_{\EX}(\hat u)
    & = \iu \lambda_{\mathrm i} \, \hat u
    \,,
    \\
    \hat \phi_{\IM}(\hat u, \theta)
    & = \left(\lambda_{\mathrm r} - \frac{\theta}{2}\lambda_{\mathrm i}^2 \right) \hat u
    \,.
  \end{align*}
%\end{subequations}
%%%
With these ingredients, the one-stage first-order method SI1(1) yields
\begin{equation*}
  %\label{eq:SI1(1):conv-diff:model}
  \hat u^{n+1} 
  = \hat u^{n}
  + \Delta t \bigl[ \hat \phi_{\EX}(\hat u^{n})
                  + \hat \phi_{\IM}(\hat u^{n+1}, \Delta t)
             \bigr]
\end{equation*}
or, equivalently,
\begin{equation*}
  %\label{eq:SI1(1):dahlquist}
  \hat u^{n+1} 
  = \hat u^{n}
  + \iu z_{\mathrm i} \hat u^{n}
  + \left( z_{\mathrm r} - \frac{1}{2} z_{\mathrm i}^2 \right) \hat u^{n+1}
  \,.
\end{equation*}
%%%
This gives the stability function
\begin{equation}
  \label{eq:stability-function:SI1(1)}
  R_{\text{SI1(1)}}(z) 
    = \frac{1 + \iu z_i}{1 - z_r + \frac{1}{2}z_i^2}
  \,.
\end{equation}
%%%
For the two-stage method SI1(2), similar considerations lead to
\begin{equation}
  \label{eq:stability-function:SI1(2)}
  R_{\text{SI1(2)}}(z) 
    = \frac{1 + \iu z_i R_{\text{SI1(1)}}(z)}{1 - z_r + \frac{1}{2}z_i^2}
  \,.
\end{equation}
%%%
The stability domains of both methods are depicted in Fig.~\ref{fig:stability-functions:si1}.
%%%
A closer inspection of the stability functions 
\eqref{eq:stability-function:SI1(1)}
and
\eqref{eq:stability-function:SI1(2)}
yields the following result:
%%%
\begin{proposition}
The time integration methods SI1(1) and SI1(2) are $L$-stable.
\end{proposition}

\begin{figure}
\includegraphics[scale=0.52]{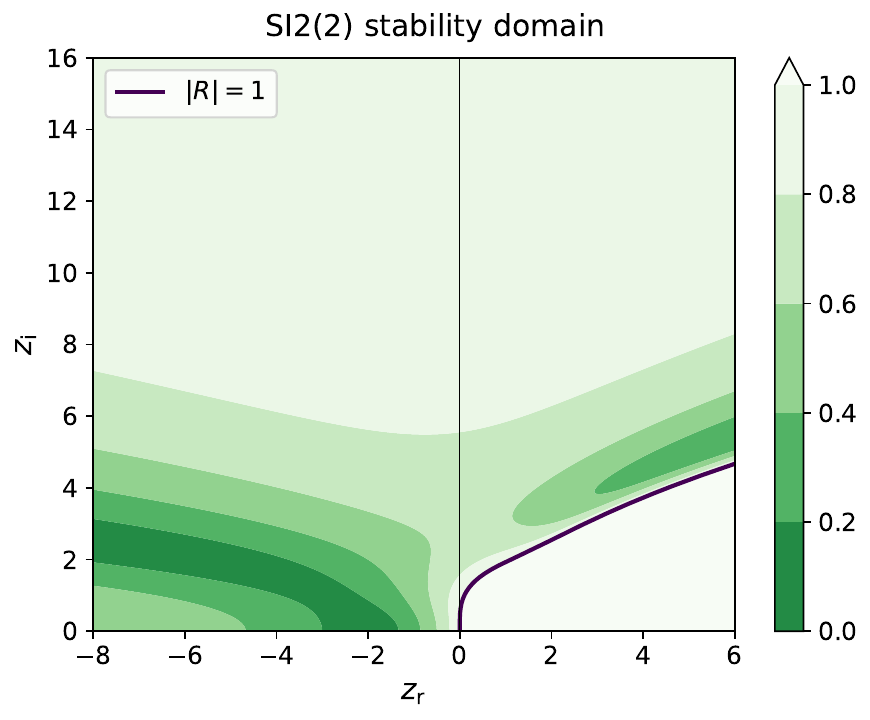}
\caption{Stability domain of the second-order integrator
\label{fig:stability-function:si22}}
\end{figure}

Finally, the stability function of the second-order method SI2(2) is given by
\begin{equation*}
  %\label{eq:stability-function:SI2(2)}
  R_{\text{SI2(2)}}(z) 
    = 1 
    + z \frac{1 + \frac{\iu}{2}z_i R_{\text{SI2(1)}}(z)}
             {1 - \frac{1}{2}z_r + \frac{1}{4}z_i^2}
  \,,
\end{equation*}
where 
\begin{equation*}
  %\label{eq:stability-function:SI2(1)}
  R_{\text{SI2(1)}}(z) 
    = \frac{1 + \frac{\iu}{2}z_i}{1 - \frac{1}{2}z_r + \frac{1}{4}z_i^2}
\end{equation*}
is the stability function of the first stage.
%%%
The graph of the stability domain depicted in Figure~\ref{fig:stability-function:si22} indicates the following property, which can be proved by analyzing the stability function:
%%%
\begin{proposition}
The time integration method SI2(1) is $A$-stable.
\end{proposition}

%===================================================================================

\pagebreak % enforced to avoid page break after subsection title

\subsection{Spectral deferred correction method}
\label{sec:time-integration:sdc}

\newcommand{\ti}{\tilde t\,}

%-----------------------------------------------------------------------------------

\subsubsection{SDC with right-sided Gauss-Radau points}
\label{sec:sdc:derivation}

The SDC method seeks an approximate solution to the Picard integral equation
\begin{equation}
  \label{eq:picard}
    \M u(t) 
      = \M u(t^n) 
      + \int_{t^n}^{t} \M f(\M u(\ti), \ti) \,\D \ti
    \,,
    \quad 
    t \in [t^n,t^{n+1}]
    \,.
\end{equation}
%%%
Its key idea is to
1) replace the exact solution and the integrand by Lagrange interpolants through a set of collocation points and
2) solve the error equation by sweeping through the collocation points with a low-order method \cite{TI_Dutt2000a}.
%%%
Common choices for the collocation points include Gauss-type quadrature points and equidistant points \cite{TI_Layton2005a,TI_Ong2020a}.
%%%
The following description is restricted to the right-sided Gauss-Radau points, but can easily be extended to other point sets.
%%%

Let ${\NM\tau = \{\tau_m\}_{m=1}^M}$ be the right Radau points in the unit interval ${[0,1]}$ and
\begin{equation}
  \label{eq:time:collocation:basis}
  \ell_{\! m}(\tau) 
  = \prod_{\substack{n=1 \\ n\ne m}}^{M} \frac{\tau - \tau_n}{\tau_m - \tau_n}
\end{equation}
the corresponding Lagrange basis polynomials.
%%%
The collocation points ${\NM t = \{t_m\}_{m=1}^M}$ in the interval ${[t^{n},t^{n+1}]}$ result from the mapping ${t_m = t(\tau_m)}$, where ${t(\tau) = t^n + \tau \Delta t}$.
%%%
The values of some quantity $\M v$ in the collocation points form the vector ${\NM{\M v} = [\M v_n]}$, which is used to construct the Lagrange interpolant
\begin{equation*}
  %\label{eq:time:collocation:interpolant}
  \M L(\NM{\M v},t) = \sum_{m = 1}^{M} \M v_m\,\ell_{\! m}(\tau(t))
\end{equation*}
with ${\tau(t) = (t - t^n)/\Delta t}$.
%%%
For a simple notation, the initial time is denoted by ${t_0 = t^n}$ and the initial value of the solution by ${\M u_0 = \M u(t^n)}$.
%%%
Note, however, that $t_0$ is no collocation point, $\M u_0$ is not a component of $\NM{\M u}$ and generally ${\M u_0 \ne \M L(\NM{\M u},t_0)}$.

The SDC method starts from an initial approximation ${\NM{\M u}_{m}^{0} = [ \M u_m^0 ]}$ that results from sweeping through all subintervals ${\{[t_{m-1},t_m]\}_{m = 1}^M}$ with a suitable time integrator.
This prediction step is followed by several correction sweeps designed to reduce the error
\begin{equation*}
  %\label{eq:sdc:error-function}
  \M\delta(t) = {\M u}(t) - \M L(\NM{\M u}^k,t)
  \,,
\end{equation*}
where ${\NM{\M u}^k}$ denotes the $k$-th approximation.
%%%
Using the Picard formulation \eqref{eq:picard} and the residual function
\begin{equation}
  \label{eq:sdc:residual-function}
  \M\varepsilon(t) 
  = \M u_0 
  + \int_{t_0}^{t}\M f\bigl(\M L(\NM{\M u}^k, \ti), \,\ti \bigr) \D \ti
  - \M L(\NM{\M u}^k,t)
\end{equation}
the error equation
\begin{equation*}
  %\label{eq:sdc:error-equation}
  \M\delta(t) 
  = \M\varepsilon(t)
  + \int_{t_0}^{t} 
      \Bigl[ \M f\bigl(\M L(\NM{\M u}^k, \ti) + \M\delta(\ti), \,\ti \bigr)
           - \M f\bigl(\M L(\NM{\M u}^k, \ti), \ti \bigr)
      \Bigr] \D \ti
\end{equation*}
can be derived, see e.g. \cite{TI_Minion2003b}.
%%%
The difference of $\M\delta(t)$ at the subinterval end points results in the identity
\begin{equation}
  \label{eq:sdc:error-equation:subinterval}
  \M\delta(t_{m}) 
  = \M\delta(t_{m-1}) 
  + \M\varepsilon(t_{m})
  - \M\varepsilon(t_{m-1})
  + \int_{t_{m-1}}^{t_{m}} 
      \Bigl[ \M f\bigl(\M L(\NM{\M u}^k, \ti) + \M\delta(\ti), \,\ti \bigr)
           - \M f\bigl(\M L(\NM{\M u}^k, \ti), \ti \bigr)
      \Bigr] \D \ti
  .
\end{equation}
%%%
Letting ${\M\delta_{m}}$ denote the approximation to $\M\delta(t_{m})$, the discretization of \eqref{eq:sdc:error-equation:subinterval} with the one-stage first-order integrator SI1(1) is 
\begin{equation}
  \begin{aligned}
    \label{eq:sdc:SI1(1):delta}
    \M\delta_{m} 
    & = \M\delta_{m-1} 
      + \M\varepsilon(t_{m})
      - \M\varepsilon(t_{m-1})
    \\
    & + \Delta t_m 
          \left[ \fex\!\left(\M u_{m-1}^{k+1}\right)
               + \fim\!\left(\M u_{m-1}^{k+1},\M u_{m}^{k+1},t_m,\Delta t_m\right)
          \right]
    \\
    & - \Delta t_m 
          \left[ \fex\!\left(\M u_{m-1}^{k}\right)
               + \fim\!\left(\M u_{m-1}^{k},\M u_{m}^{k},t_m,\Delta t_m\right)
          \right]
    \,,
  \end{aligned}
\end{equation}
where ${\Delta t_m = t_m - t_{m-1}}$.
%%%
With \eqref{eq:sdc:residual-function}, the residual difference can be written
\begin{equation}
  \label{eq:sdc:residual:difference}
  \M\varepsilon(t_{m}) - \M\varepsilon(t_{m-1})
  = \int_{t_{m-1}}^{t_m}\!\M f\bigl(\M L(\NM{\M u}^k,\ti),\ti\bigr) \D \ti
  - \M u_{m}^k
  + \M u_{m-1}^k
  \,.
\end{equation}
The integral in \eqref{eq:sdc:residual:difference} is approximated using the quadrature
\begin{equation*}
  %\label{eq:sdc:quadrature}
  \int_{t_{m-1}}^{t_m}\!\M f\bigl(\M L(\NM{\M u}^k,t),t\bigr) \D t
  \approx
    \int_{t_{m-1}}^{t_m}\!\M L\left(\NM{\M f}^k,t\right) \D t
  = \Delta t \sum_{i=1}^{M} w_{i,m} \M f(\M u_i^k,t_i)
  = \M S_m^k
\end{equation*}
with the weights
\begin{equation*}
  %\label{eq:sdc:weights}
  w_{i,m} = \int_{0}^{\tau_m}\! \ell_{\! i}(\tau) \D \tau
  \,.
\end{equation*}
Substituting \eqref{eq:sdc:residual:difference} into \eqref{eq:sdc:SI1(1):delta} and defining the new approximation ${\M u_m^{k+1} = \M u_m^{k} + \M\delta_{m}}$ gives the SI1(1) corrector
\begin{equation}
  \begin{aligned}
    \label{eq:sdc:corrector:SI1(1):semi}
    \M u_{m}^{k+1} 
    & = \M u_{m-1}^{k+1} 
      + \M S_m^k
    \\
    & + \Delta t_m 
          \left[ \fex\!\left(\M u_{m-1}^{k+1}\right)
               + \fim\!\left(\M u_{m-1}^{k+1},\M u_{m}^{k+1},t_m,\Delta t_m\right)
          \right]
    \\
    & - \Delta t_m 
          \left[ \fex\!\left(\M u_{m-1}^{k}\right)
               + \fim\!\left(\M u_{m-1}^{k},\M u_{m}^{k},t_m,\Delta t_m\right)
          \right]
  \end{aligned}
\end{equation}
for 
${m =1, \dots M}$.

Similarly, applying the two-stage first-order integrator \eqref{eq:SI1(2):semi} to \eqref{eq:sdc:error-equation:subinterval} results in the SI1(2) corrector
\begin{equation*}
  %\label{eq:sdc:corrector:SI1(2):semi}
  \begin{alignedat}{4}
    &  \M u^{(1)} 
    &&  = \M u_{m-1}^{k+1} 
        + \M S_m^k
    &&  + \Delta t_m
            \bigl[ \fex\!\left(\M u_{m-1}^{k+1}\right)
    &&           + \fim\!\left(\M u_{m-1}^{k+1},\M u^{(1)},t_m,\Delta t_m\right)
            \bigr]
    \\
    &
    &&
    &&  - \Delta t_m
             \bigl[ \fex\!\left(\M u_{m-1}^{k}\right)
    &&            + \fim\!\left(\M u_{m-1}^{k},\M u_{m}^{k},t_m,\Delta t_m\right)
             \bigr]      
    \,,
    \\
    &  \M u^{(2)} 
    &&  = \M u_{m-1}^{k+1} 
        + \M S_m^k
    &&  + \Delta t_m
            \bigl[ \fex\!\left(\M u^{(1)}\right)
    &&           + \fim\!\left(\M u_{m-1}^{k+1},\M u^{(2)},t_m,\Delta t_m\right)
            \bigr]
    \\
    &
    &&
    &&  - \Delta t_m
             \bigl[ \fex\!\left(\M u_{m}^{k}\right)
    &&            + \fim\!\left(\M u_{m-1}^{k},\M u_{m}^{k},t_m,\Delta t_m\right)
             \bigr]      
    \,,
    \\
    &  \M u_{m}^{k+1} 
    && = \M u^{(2)} 
  \end{alignedat}
\end{equation*}
for 
${m =1, \dots M}$.

\begin{remark}
Both correctors and thus the resulting SDC methods are semi-implicit, as the intermediate and final solutions appear only on the left-hand side and in $\fim$. 
\end{remark}

It should be noted that a similar corrector can be derived from SI2(2).
%%%
However, this approach would require additional data in the subinterval midpoints and would not promise faster convergence with non-equidistant collocation points \cite{TI_Christlieb2009a}.
%%%
Therefore, the present study focuses on SDC methods using first-order predictors and correctors.
%%%
For concise notation, these methods are labeled SDC-SI($s_1$,$s_2$)$_M^K$, where
$s_1$ specifies the number of predictor stages, 
$s_2$ the number of corrector stages,
$M$   the number of collocation points and
$K$   the total number of iterations, 
i.e. one predictor sweep and ${K-1}$ corrector sweeps.
%%%
\begin{remark}
When applying the integrators SI1(1) or SI1(2) as a predictor, they sweep through all subintervals one after the other. As a result the step size $\Delta t$ is replaced by the subinterval length ${\Delta t_m}$ in \eqref{eq:SI1(1):semi} and
\eqref{eq:SI1(2):semi}, respectively. Note that $\theta$ simultaneously assumes the same value, as it is fixed to the step size.
\end{remark}

%-----------------------------------------------------------------------------------

\subsubsection{Equivalence to Radau IIA collocation and discontinuous Galerkin methods}
\label{sec:sdc:collocation+dg}

If the SDC method converges, the difference between $\M u_{m}^{k}$, $\M u_{m}^{k+1}$ and $\M u_{m}$ vanishes for ${k \rightarrow \infty}$.
%%%
As a consequence, the solution satisfies the collocation scheme
\begin{equation}
  \label{eq:collocation:incremental}
  \M u_{m} = \M u_{m-1} + \Delta t \sum_{i=1}^{M} w_{i,m} \M f(\M u_i,t_i)
  \,,\quad
  m = 1,\dots M
  \,.
\end{equation}
By transforming \eqref{eq:collocation:incremental} into a non-incremental form it can be shown that the scheme is identical to the Radau IIA collocation method of the order ${2M-1}$
\begin{equation}
  \label{eq:collocation:radau-iia:semi}
  \M u_{m} = \M u(t_0) 
           + \Delta t \sum_{n=1}^{M} a_{m,n} \M f(\M u_n,t_0 + c_n \Delta t)
  \,,\quad
  m = 1,\dots M
\end{equation}
with the nodes ${c_n = \tau_n}$ and coefficients
\begin{equation*}
  a_{m,n}
  = \sum_{k=1}^{m} w_{n,k}
  = \int_{0}^{\tau_m}\! \ell_{\! n}(\tau) \D \tau
  \,,
\end{equation*}
see \citet[Ch.~6]{TI_Deuflhard2002a}.
%%%
It is therefore expected that the corresponding SDC method will achieve the same order.
%%%
Furthermore, it was proven in \cite{TI_Huynh2023a} that the Radau IIa method and thus the SDC method is equivalent to a discontinuous Galerkin method using the basis functions \eqref{eq:time:collocation:basis}.
% together with an upwind flux to ensure temporal causality

%-----------------------------------------------------------------------------------

\subsubsection{Extension to other point sets and relation to Euler-based SDC methods}
\label{sec:sdc:extension+relations}

Similar SDC methods can be constructed using other sets of collocation points, e.g.,  Gauss, Gauss-Lobatto or equidistant points.
%%%
The latter two include the left interval boundary, which must be taken into account when constructing the Lagrange interpolant $\M L(\NM{\M u},t)$.
%\eqref{eq:time:collocation:interpolant}
%%%
In addition, the first coefficient is then determined by the initial condition.
%%%
Furthermore, it is noted that the one-step corrector \eqref{eq:sdc:corrector:SI1(1):semi} reproduces the Euler-based semi-implicit SDC method proposed by \citet{TI_Minion2003b} when ${\theta = 0}$ is enforced in $\fim$.
%%%
This method was shown to increase the order of the approximate solution by one with each correction sweep, until the order of the underlying collocation method is reached \cite{TI_Hagstrom2006a,TI_Hansen2011a,TI_Causley2019a}.
%%%
In the following, the method is referred to as SDC-EU and used for comparison.
%%%
%Due to the similarity, it is assumed that the proposed SDC-SI methods converge at the same %rate as SDC-EU.
Due to the similarity between SDC-EU and SDC-SI, it is assumed that both methods converge at the same rate.

%-----------------------------------------------------------------------------------

\subsubsection{Linear stability of SDC methods}
\label{sec:sdc:stability}

%-----------------------------------------------------------------------------------

The Gauss-type collocation methods are known to be $A$-stable and the Radau IIA method is even $L$-stable \cite{TI_Hairer1996a}.
%%%
Unfortunately, corresponding SDC methods do not necessarily inherit these properties.
%%%
Even using backward Euler for prediction and correction is not sufficient to achieve $A$-stability \cite{TI_Dutt2000a}.
%%%
Nevertheless, stable semi-implicit SDC methods can be constructed using the correctors proposed above.
%%%

\begin{figure}
  \subcaptionbox{SDC with IMEX-Euler predictor and corrector
    \label{fig:stability-domains:sdc:euler}}
    {\includegraphics[scale=0.52]{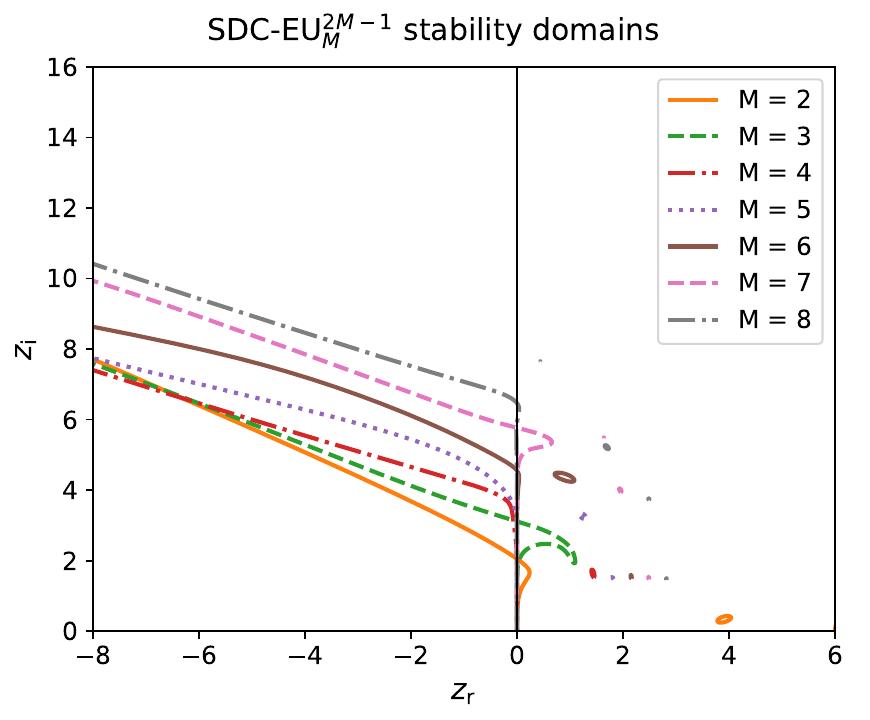}}
  \subcaptionbox{SDC with SI1(1) predictor and corrector
    \label{fig:stability-domains:sdc:si11}}
    {\includegraphics[scale=0.52]{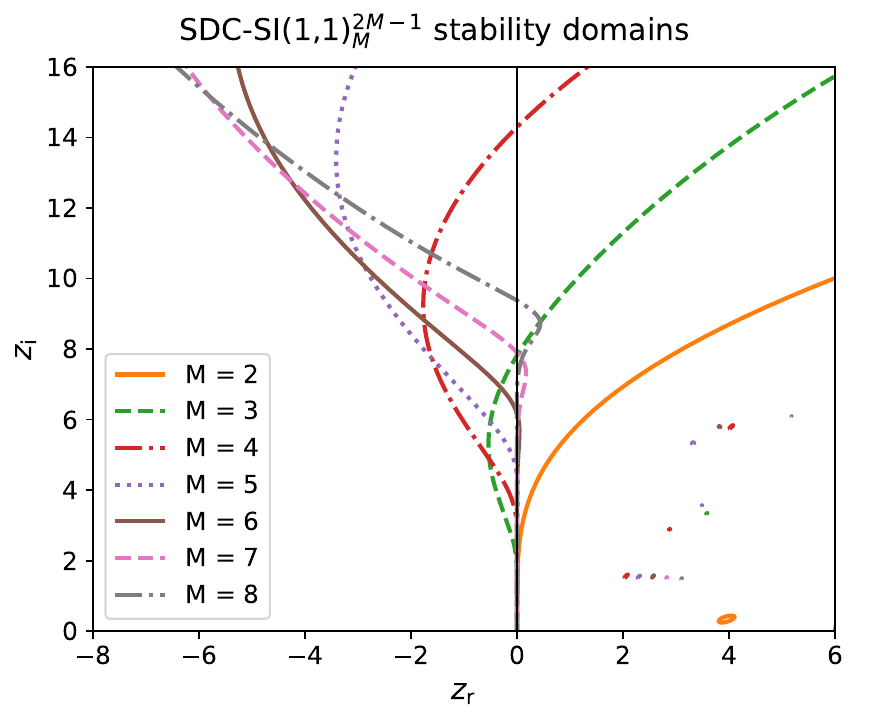}} 
  \\[\medskipamount]
  \subcaptionbox{SDC with SI1(2) predictor and corrector
    \label{fig:stability-domains:sdc:si12}}
    {\includegraphics[scale=0.52]{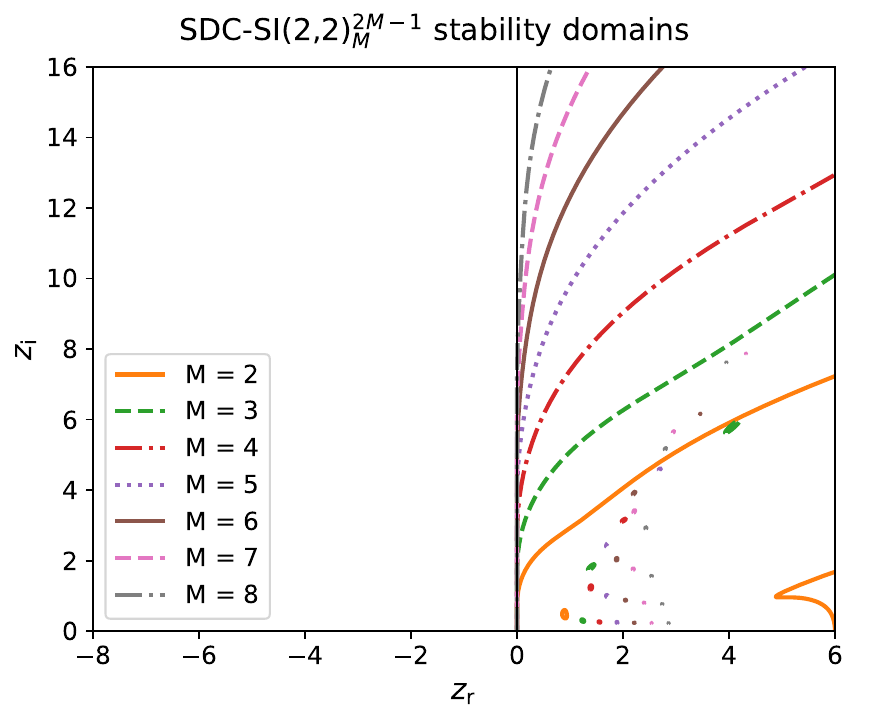}}
  \subcaptionbox{SDC with SI1(2) close to imaginary axis
    \label{fig:stability-domains:sdc:si12:zoom}}
    {\includegraphics[scale=0.52]{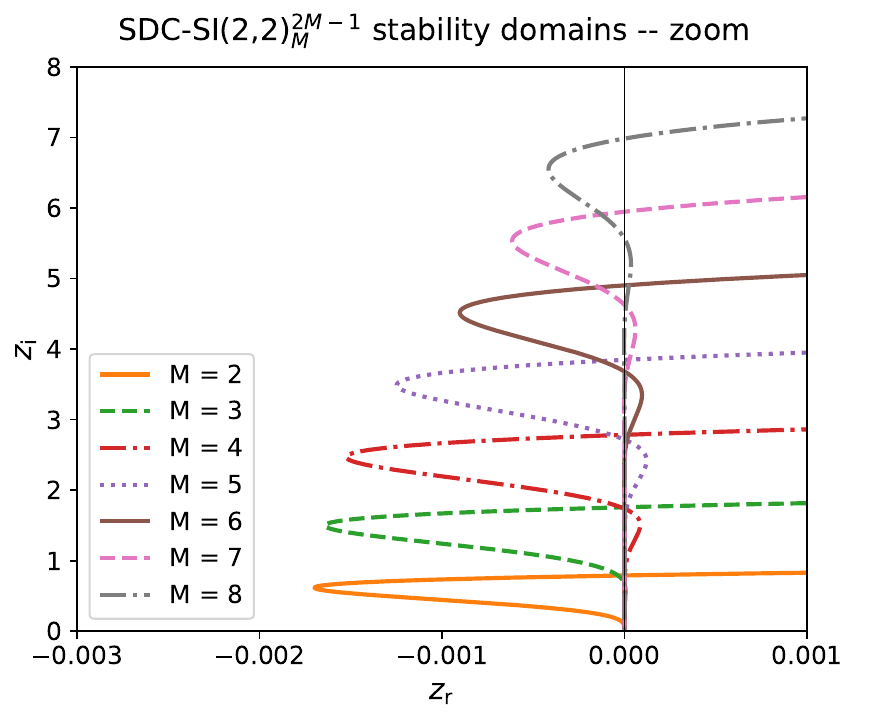}}
  \\[\medskipamount]
  \subcaptionbox{SDC optimized for stability
    \label{fig:stability-domains:sdc:si1s}}
    {\includegraphics[scale=0.52]{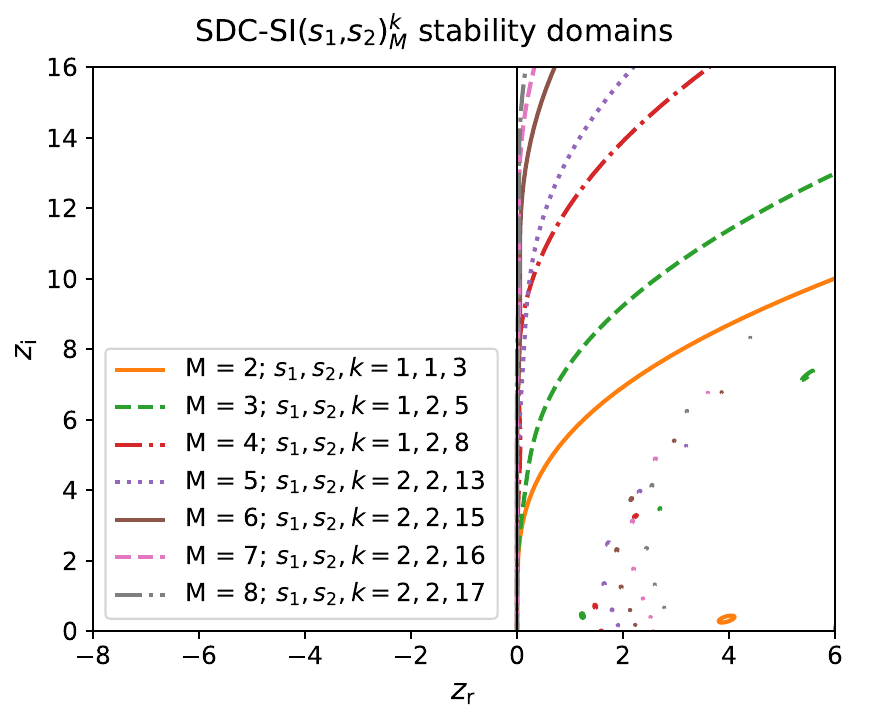}} 
  \subcaptionbox{Optimized SDC close to imaginary axis
    \label{fig:stability-domains:sdc:si1s:zoom}}
    {\includegraphics[scale=0.52]{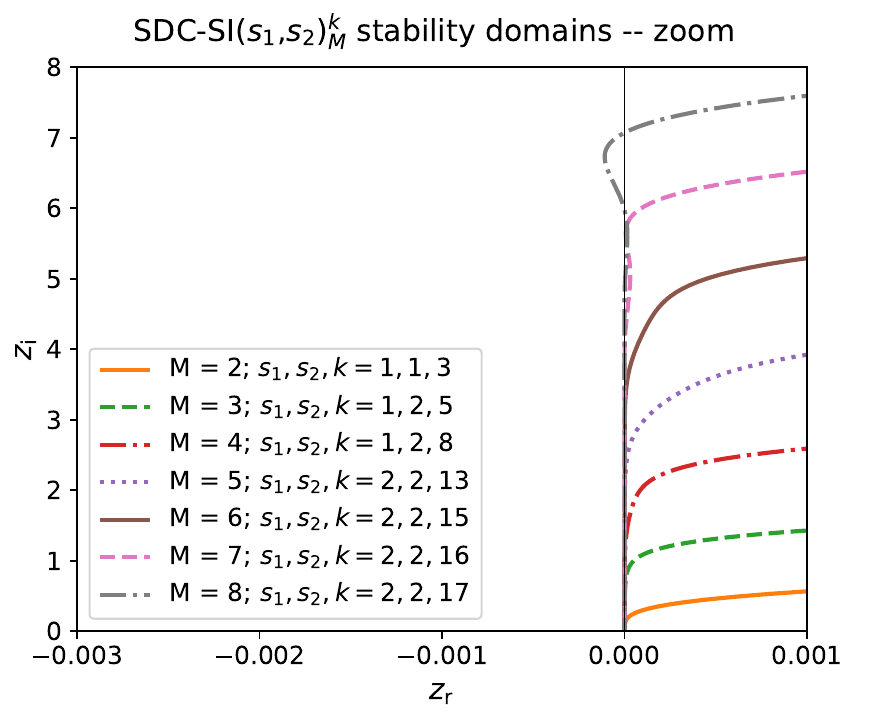}} 
  \caption{Neutral stability curves for SDC methods with up to
    ${M\!=\!8}$ Radau points using different predictors and correctors
    %, ${|R(z)| = 1}$,
    \label{fig:stability-domains:sdc}}
\end{figure}

Although the stability function of an SDC method is difficult to determine explicitly, it can be evaluated numerically by executing one time step of size ${\Delta t = 1}$ for the Dahlquist equation \eqref{eq:dahlquist} with the initial condition ${\hat u(0) = 1}$.
%%%
For comparison, Figure~\ref{fig:stability-domains:sdc:euler} shows the stability domains of the Euler-based semi-implicit SDC methods using ${M = 2}$ to $8$ Radau points and ${K = 2M - 1}$ iterations.
%%%
The slope of the neutral stability curves suggests that these methods are $A(\alpha)$-stable with $\alpha$ between 18° and 40°.
%%%
In the case of pure convection, stability requires ${z_{\mathrm i} \lesssim M/2}$, which is disappointing when compared with IMEX Runge-Kutta methods \cite{TI_Ascher1997a,TI_Cavaglieri2015a}.
%%%
Using the SI1(1) integrator for prediction and correction improves the stability considerably (Fig.~\ref{fig:stability-domains:sdc:si11}).
%%%
The third-order method with two Radau points is even $L$-stable within the scope of numerical accuracy.
%%%
This property is lost with increasing $M$, however 
the methods remain $A(\alpha)$-stable with $\alpha$ between 40° and 80° and 
the convective stability limit doubles to ${z_{\mathrm i} \lesssim M}$.
%%%
The SDC methods, which are based on the two-stage integrator SI1(2), achieve a further improvement and even appear $A$-stable at first glance (Fig.~\ref{fig:stability-domains:sdc:si12}).
%%%
However, the zoom in Fig.~\ref{fig:stability-domains:sdc:si12:zoom} reveals unstable regions near the imaginary axis.
%%%
A further improvement of stability was achieved by combining different predictors and correctors and adjusting the number of iterations.
%%%
The resulting parameters are given in Tab.~\ref{tab:sdc:stable:overview}.
%%%
Figure~\ref{fig:stability-domains:sdc:si1s} shows the corresponding stability domains and Fig.~\ref{fig:stability-domains:sdc:si1s:zoom} a close-up near the imaginary axis.
%%%
Detailed numerical studies suggest that the improved methods are $L$-stable up to order 11 or 6 Radau points.
%%%
The methods with 7 and 8 points are not $A$-stable, but instability is observed only in small regions close to the imaginary axis, see Tab.~\ref{tab:sdc:stable:overview}.

\begin{table}
  \begin{tabular}{@{}cccccc@{}} \toprule
  order &  $M$  & $s_1$ & $s_2$ &  $k$  & $z_{r,\max}$     \\ \midrule
    3   &   2   &   1   &   1   &   3   &  0               \\       
    5   &   3   &   1   &   2   &   5   &  0               \\       
    7   &   4   &   1   &   2   &   8   &  0               \\       
    9   &   5   &   2   &   2   &  13   &  0               \\       
   11   &   6   &   2   &   2   &  15   &  0               \\       
   13   &   7   &   2   &   2   &  16   &  \textminus5.2e\textminus7 \\
   15   &   8   &   2   &   2   &  17   &  \textminus1.1e\textminus4 \\ \bottomrule
  \end{tabular}
\caption{SDC-SI($s_1$,$s_2$)$_M^k$ methods with optimal stability properties.
  The first column indicates the order of the underlying Gauss-Radau IIA method,
  $M$ is the number of collocation points,
  $s_1$ the number of predictor stages,
  $s_2$ the number of corrector stages,
  $k$ the total number of iterations and
  $z_{r,\max}$ the largest value on the real axis, so that the method is stable
  for all $z_i$.
  \label{tab:sdc:stable:overview}}
\end{table}

%===================================================================================

% !TEX encoding = UTF-8 Unicode
% !TEX root = ssi-sdc-1d.tex

%%%%%%%%%%%%%%%%%%%%%%%%%%%%%%%%%%%%%%%%%%%%%%%%%%%%%%%%%%%%%%%%%%%%%%%%%%%%%%%%%%%%

\section{Spatial discretization}
\label{sec:spatial-discretization}

%===================================================================================

\subsection{Preliminaries}
\label{sec:spatial-discretization:preliminaries}

The time integration methods developed in the previous section can be combined with any spatial discretization method.
%%%
In the following, the discontinuous Galerkin method with nodal basis functions is used \cite{SE_Hesthaven2008a}.
%%%
Subdividing $\Omega$ into non-overlapping elements $\Omega^e$ yields the computational domain ${\Omega_h = \cup\Omega^e}$ and the related set of element boundaries $\Gamma_h$.
%%%
For a given time $t$, the approximate solution $\M u_h(x,t)$ is sought in the function space
\begin{equation*}
  \mathbb U_h = \left\{
    \M u_h \in \bigl[\mathbb L^2(\Omega)\bigr]^d  : 
    \left. \M u_h \right|_{\Omega^e} 
    \in  \bigl[\mathbb P_{\! P}(\Omega^e)\bigr]^d 
    \quad 
    \forall \Omega^e \in \Omega
  \right\},
\end{equation*}
where
$\mathbb L^2$ is the space of square-integrable functions and 
$\mathbb P_{\! P}$ the space of polynomials with degree less than or equal to $P$.
%%%
The test space $\mathbb W_h$ is mathematically identical to $\mathbb U_h$, but has reciprocal physical dimensions, so that the product ${\transpose{\M w_h}\M u_h}$ is well-defined for any ${\M w_h \in \mathbb W_h}$.
%%%
Finally, the average and jump operators on the joint boundary between elements $\Omega^-$ and $\Omega^+$ are defined by
\begin{alignat*}{2}
  & \avg{\M u_h} && = \frac{\M u_h^- + \M u_h^+}{2} \,, \\
  & \jmp{\M u_h} && = n^- \M u_h^- + n^+ \M u_h^+   \,,
\end{alignat*}
where 
$\M u_h^{\pm}$ are the traces of the corresponding element solutions and
$n^{\pm}$ the normals.

%-----------------------------------------------------------------------------------

%===================================================================================

\subsection{Discontinuous Galerkin formulation}
\label{sec:spatial-discretization:dg-formulation}

Terms with no spatial derivatives are discretized using the projection operator
\begin{equation}
  \label{eq:space:dg:projection}
  \PO(\M w_h,\M u_h) = \int_{\Omega_h} \transpose{\M w_h} \M u_h\,\D x
  \,.
\end{equation}
%%%
In particular, the source term becomes
\begin{equation*}
  %\label{eq:space:dg:sources}
  \Fs(\M w_h,\M u_h, t) = \PO(\M w_h,\fs(\M u_h,x,t))
  \,.
\end{equation*}
%%%
The convection term is discretized as
\begin{equation}
  \label{eq:space:dg:convection}
  \Fc(\M w_h, \M u_h) 
    = \int_{\Omega_h} \transpose{(\d_x \M w_h)} \fc(\M u_h)\,\D x
    + \sum_{\Gamma_h} \transpose{\jmp{\M w_h}}\widehat{\fc}(\M u_h)
    \,,
\end{equation}
where $\widehat{\fc}$ is the Riemann flux in the scalar case and the Roe flux with sonic fix of Harten and Hyman \cite{FV_Harten1983a} in the Euler or Navier-Stokes case.
%%%
Applying the symmetric interior penalty method to diffusion terms results in
%%%
\begin{equation}
  \label{eq:space:dg:diffusion}
  \begin{aligned}
    \Fd(\M w_h,\M u_h, \M A_h)
      = &- \int_{\Omega_h} \transpose{\d_x \M w_h} \M A_h \d_x \M u_h  \,\D x
     \\ &+ \sum_{\Gamma_h} 
           \bigl( \avg{\transpose{\M A_h} \d_x \M w_h} \jmp{\M u_h}
                + \jmp{\M w_h} \avg{\M A_h \d_x \M u_h}
           \bigr) \\
        & - \sum_{\Gamma_h} \mu\jmp{\M w_h}\M{{\widehat{A}}}\jmp{\M u_h}
     \,,
  \end{aligned}
\end{equation}
where
${\M A_h}$ is the diffusion matrix, $\M{{\widehat{A}}}$ the interface diffusion matrix with components
\begin{equation*}
  %\label{eq:space:dg:diffusion:interface}
  \widehat A_{i,j} = 
  \begin{cases}
    \max(A_{i,j}^-,A_{i,j}^+)           & \quad i = j \\
    \frac{1}{2}(A_{i,j}^- + A_{i,j}^+)  & \quad i \ne j
  \end{cases}
\end{equation*}
and
\begin{equation}
  \label{eq:space:dg:diffusion:penalty}
  \mu = c_{\mu} \avg{\frac{P(P+1)}{2\Delta x^e}}
\end{equation}
the penalty factor with the element width ${\Delta x^e}$ and ${c_{\mu} > 1}$, see e.g.
\cite{SE_Stiller2016b}.

\begin{remark}
\label{rem:conservation}
The numerical fluxes introduced in the convection and diffusion terms are unique at the boundaries between the elements and thus preserve the conservation properties.
%%%
This can illustrated by choosing ${\M w_h = \MC 1}$ which eliminates all boundary contributions in \eqref{eq:space:dg:convection} and \eqref{eq:space:dg:diffusion}.
\end{remark}

%===================================================================================

\subsection{Shock capturing}
\label{sec:spatial-discretization:shock-capturing}

\newcommand{\vs} {\nu_{\mathrm s}}
\newcommand{\cs} {c_{\mathrm s}}
\newcommand{\ds} {\kappa_{\mathrm s}}

Following \citet{SE_Persson2006a}, an artificial diffusivity $\vs$ is introduced to capture unresolved discontinuities.
%%%
The diffusivity is constant in each element and applies to all conservation variables.
%%%
It can be included by adding a diagonal entry to the diffusion matrix, so that  ${\fd = (\Ad + \vs \MC I) \,\d_x \M u}$.
%%%
The method is briefly described to define the parameters: 
%%%

First select a quantity of interest $q$ and determine its Legendre expansion
\begin{equation*}
  q^e = \sum_{i=0}^{P} q_i^e L_i
\end{equation*}
in element $\Omega^e$.
Define the filtered quantity
\begin{equation*}
  \tilde q^e = \sum_{i=0}^{P-1} q_i^e L_i
\end{equation*}
and evaluate the smoothness indicator
\begin{equation*}
  s^e = \log_{10}\left(\frac{\norm{q^e - \tilde q^e}_2^2}{\norm{q^e}_2^2}\right),
\end{equation*}
where $\norm{\cdot}$ is the $L^2$ norm on $\Omega^e$.
%%%
Then define the artificial diffusivity in $\Omega^e$ by
\begin{equation}
  \label{eq:art-diffusivity}
  \nu_{\mathrm s}^e 
  = \begin{cases}
       0     & \text{if}  \quad  s_0 - \ds > s^e \,, \\
       \hat\nu_{\mathrm s}^e 
              \,\frac{1}{2}\left(1 + \sin\frac{\pi(s^e - s_0)}{2\ds}\right)
             & \text{if}  \quad  s_0 - \ds \le s^e \le s_0 + \ds \,, \\
       \hat\nu_{\mathrm s}^e  
             & \text{if} \quad
                 s_0 + \ds < s^e \,,
    \end{cases}
\end{equation}
where
${s_0}$ is a reference value,
$\ds$ is the half width of the activation interval,
\begin{equation*}
  %\label{eq:nu_s:max}
  \hat\nu_{\mathrm s}^e = \cs \lambda_{\mathrm c, \max}^e \frac{\Delta x^e}{P}
\end{equation*}
the maximum diffusivity, $\cs$ a scaling factor and ${\lambda_{\mathrm c, \max}^e}$ the maximum magnitude of the eigenvalues of $\Ac$ in $\Omega^e$.
%%%
In the studies presented below, ${q = u}$ was chosen for the Burgers equation and ${q = \rho}$ for the Euler equations.
%%%
The reference value is fixed to ${s_0 = -4(\log_{10} P + 1)}$, while $\cs$ and $\ds$ are adapted to each test case.
%%%
It should be noted that several improvements to the method have been proposed, but these are beyond the scope of this paper, see e.g. \cite{SE_Barter2010a,SE_Kloeckner2011a,SE_Glaubitz2019a}.

%===================================================================================

\subsection{Application to time integration methods}
\label{sec:spatial-discretization:application}

In the course of spatial discretization, the RHS functions introduced in Sec.~\ref{sec:time-integration} are replaced by the functionals
\begin{align*}
  \Fex\bigl(\M w_h, \M u_h^\alpha\bigr) 
  &= \Fc\bigl(\M w_h, \M u_h^\alpha\bigr)
  \,,
  \\
  \Fim\bigl(\M w_h, \M u_h^\alpha, \M u_h^\beta, t, \theta\bigr) 
  &= \Fd\bigl( \M w_h, \M u_h^\beta, 
               \tfrac{\theta}{2}\Ac^2(\M u_h^\alpha)
               + \Ad(\M u_h^\alpha) 
               + \nu_s(\M u_h^\alpha) \MC I
        \bigr)
   + \Fs\bigl(\M w_h,\M u_h^\alpha, t\bigr)
  \,,
  \\
  \F\bigl(\M w_h, \M u_h^\alpha, t \bigr) 
  &= \Fc\bigl(\M w_h, \M u_h^\alpha\bigr)
   + \Fd\bigl( \M w_h, \M u_h^\alpha, \Ad(\M u_h^\alpha) + \nu_s(\M u_h^\alpha) \MC I \bigr)
   + \Fs\bigl(\M w_h,\M u_h^\alpha, t\bigr)
  \,.
\end{align*}
%%%
As a result, the SI1(1) integrator \eqref{eq:SI1(1):semi} takes form
\begin{equation*}
  %\label{eq:SI1(1):full}
  \PO(\M w_h,\M u_h^{n+1})
  = \PO(\M w_h,\M u_h^{n})
  + \Delta t \bigl[ \Fex(\M w_h, \M u_h^{n})
                  + \Fim(\M w_h, \M u_h^n, \M u_h^{n+1}, t^{n+1}, \Delta t)
             \bigr]
\end{equation*}
and the SI1(1) corrector \eqref{eq:sdc:corrector:SI1(1):semi} becomes
\begin{equation*}
  \begin{aligned}
    %\label{eq:sdc:corrector:SI1(1):full}
    \PO(\M w_h,\M u_{h,m}^{k+1})
    & = \PO(\M w_h,\M u_{h,m-1}^{k+1}) 
      + \mathcal S_{h,m}^k
    \\
    & + \Delta t_m 
          \bigl[ \Fex\bigl(\M w_h, \M u_{h,m-1}^{k+1}\bigr)
               + \Fim\bigl(\M w_h, \M u_{h,m-1}^{k+1},\M u_{h,m}^{k+1},t_m,\Delta t_m\bigr)
          \bigr]
    \\
    & - \Delta t_m 
          \bigl[ \Fex\bigl(\M w_h, \M u_{h,m-1}^{k}\bigr)
               + \Fim\bigl(\M w_h, \M u_{h,m-1}^{k},\M u_{h,m}^{k},t_m,\Delta t_m\bigr)
          \bigr]
  \end{aligned}
\end{equation*}
with
\begin{equation*}
  \mathcal S_{h,m}^k 
  = \PO(\M w_h, \M S_{m}^k)
  = \Delta t \sum_{i=1}^{M} w_{i,m} \F(\M w_h, \M u_{h,i}^{k}, t_i)
  \,.
\end{equation*}
%%%
The application to the two-stage methods proceeds analogously.
%%%
Finally, the discrete form of the collocation method \eqref{eq:collocation:radau-iia:semi} reads
\begin{equation}
  \label{eq:collocation:radau-iia:full}
  \PO(\M w_h,\M u_{h,m}) = \PO(\M w_h,\M u_h(t_0))
      + \Delta t \sum_{n=1}^{M} a_{m,n} \F(\M w_h, \M u_{h,n}, t_0 + c_n \Delta t),
  \quad
  m = 1,\dots M
  \,.
\end{equation}
\begin{remark}
Equation \eqref{eq:collocation:radau-iia:full} corresponds to a space-time discontinuous Galerkin formulation due to the equivalence between Radau IIA and DG in time.
\end{remark}
\begin{remark}
The collocation method and the SDC methods presented above are conservative, as they are linear combinations of the conservative building blocks $\Fex$ and $\Fim$.
\end{remark}

%%%

%===================================================================================

% !TEX encoding = UTF-8 Unicode
% !TEX root = ssi-sdc-1d.tex

%%%%%%%%%%%%%%%%%%%%%%%%%%%%%%%%%%%%%%%%%%%%%%%%%%%%%%%%%%%%%%%%%%%%%%%%%%%%%%%%%%%%

\section{Solution techniques}
\label{sec:solution}

%===================================================================================

\subsection{Numerical integration}
\label{sec:solution:quadrature}

The integrals in the spatial operators are evaluated using a Gauss-Lobatto quadrature
with $Q+1$ points.
%%%
For the convection term $Q$ is chosen large enough to avoid aliasing errors caused by nonlinearities.
%%%
In the numerical experiments reported below,
${Q=P}$ is used for the convection-diffusion equation,
${Q=\lceil 3P/2 \rceil}$ for the Burgers equation and
${Q=2P}$ for the Euler and Navier-Stokes equations.
%%%
The remaining terms are integrated with ${Q=P}$ so that the quadrature points coincide with the spatial collocation points and therefore no interpolation is necessary.
%%%
This approach leads to a diagonal mass matrix, which simplifies the evaluation and inversion of the projection operator \eqref{eq:space:dg:projection}.
%%%

%===================================================================================

\subsection{Solution of implicit equations}
\label{sec:solution:implicit-equations}

The spatio-temporal discretization leads to linear systems of equations that can be written symbolically in the form
\begin{equation*}
  %\label{eq:solution:linear-system}
  \NM{\M K} \NM{\M u} = \NM{\M b}
  \,,
\end{equation*}
where 
$\NM{\M u}$ is the solution vector,
$\NM{\M K}$ a positive definite matrix and
$\NM{\M b}$ the RHS.
%%%
Systems of this type can be solved very efficiently with tailored multigrid methods
\cite{SE_Lottes2005a,SE_Janssen2011a,SE_Stiller2016a,SE_Stiller2016b}, but these are beyond the scope of the present study.
%%%
Instead, a preconditioned conjugate gradient method is used for scalar cases and the FGMRES method in the Euler and Navier-Stokes cases \cite{KR_Hestenes1952a,KR_Shewchuk1994a,KR_Golub1999a}.
%%%
Additionally, a fast direct solver was developed for the convection diffusion equation with constant diffusion matrices ${\Ad}$ and $\Ac$.
%%%
This solver is used in performance tests to demonstrate the competitiveness of the proposed SDC methods.
%%%

%===================================================================================

% !TEX encoding = UTF-8 Unicode
% !TEX root =  ssi-sdc-1d.tex

%%%%%%%%%%%%%%%%%%%%%%%%%%%%%%%%%%%%%%%%%%%%%%%%%%%%%%%%%%%%%%%%%%%%%%%%%%%%%%%%%%%%

\section{Numerical Experiments}
\label{sec:numerical-experiments}

%===================================================================================

The accuracy and stability of the newly developed SDC-SI methods are evaluated by means of numerical experiments for linear convection-diffusion, Burgers and compressible flow problems.
%%%
For comparison with the linear stability analysis in Sec.~\ref{sec:time-integration}, the CFL number is defined as 
\begin{equation*}
  \CFL = \frac{\Delta t\,\lambda_{\mathrm c, \max}}{\Delta x}
  \,,
\end{equation*}
where $\Delta x$ is a length scale characterizing the mesh spacing.
%%%
Following the analysis in \cite[Sec. 7.3.3]{SE_Canuto2011a}, the length scale is chosen as
\begin{equation*}
  \Delta x = \frac{\Delta x^e}{2 \delta} 
  \,,
\end{equation*}
where ${\delta \sim P^2}$ is the largest eigenvalue of the convection problem with unit velocity and one-sided Dirichlet conditions in the standard element ${[-1,1]}$.
%%%
Unless stated otherwise, the CFL number is determined from the initial conditions.
%%%
The IMEX Euler-based SDC-EU method, 
the second-order IMEX-RK CB2 Runge-Kutta method of \citet{TI_Cavaglieri2015a},
the third order explicit TVD-RK3 method of \citet{TI_Shu1988b} and 
the third order ARS(4,4,3) method of \citet{TI_Ascher1997a} 
are used for comparison.

%===================================================================================

\subsection{Convection-diffusion equation}
\label{sec:numerical-experiments:conv-diff}

%-----------------------------------------------------------------------------------

The first set of test problems relates to the convection-diffusion equation \eqref{eq:conv-diff} with the exact solution given by
\begin{equation}
  \label{eq:wave-packet}
  u(x,t) = \sum_{i=1}^7 a_i \sin\big( \kappa_i(x - s_i - vt) \big)
                        e^{-\kappa_i^2 \nu t}
  \,
\end{equation}
and ${v = 1}$.
%%%
Equation~\eqref{eq:wave-packet} defines a wave packet including a variety of wave numbers $\kappa_i$, amplitudes $a_i$ and phase shifts $s_i$ that are specified in Tab.~\ref{tab:wave-packet}.
%%%
The tests are performed in the periodic domain ${\Omega = [0,1]}$ using 64 elements of degree 15.
%%%
This resolution is fine enough to neglect spatial discretization errors.

\begin{figure}
  \subcaptionbox{error of standalone time integrators for ${\nu=0}$
    \label{fig:conv-diff:wave:nu=0e-0:base:error}}
    {\includegraphics[scale=0.52]{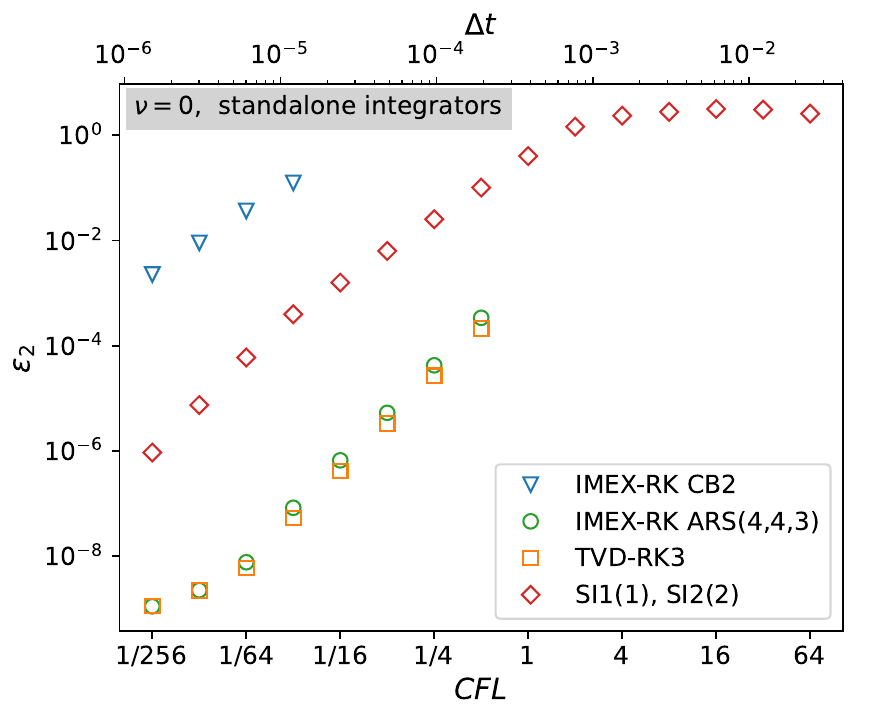}}
  \hfill
  \subcaptionbox{SDC-EU error for ${\nu=0}$
    \label{fig:conv-diff:wave:nu=0e-0:sdc:eu:error}}
    {\includegraphics[scale=0.52]{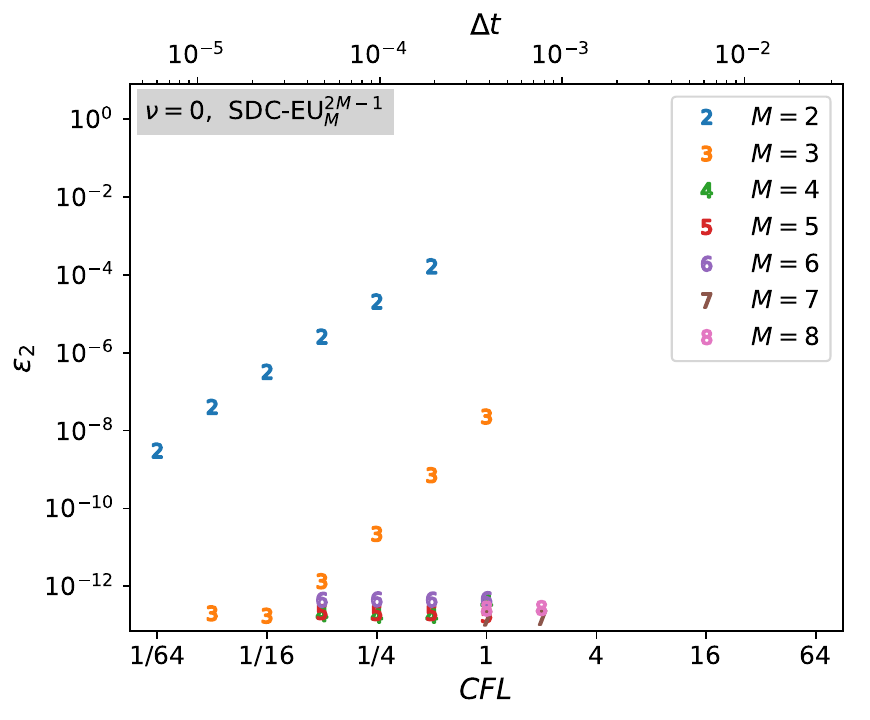}} 
  \\[\medskipamount]
  \subcaptionbox{SDC-SI error for ${\nu=0}$
    \label{fig:conv-diff:wave:nu=0e-0:sdc:si1s:error}}
    {\includegraphics[scale=0.52]{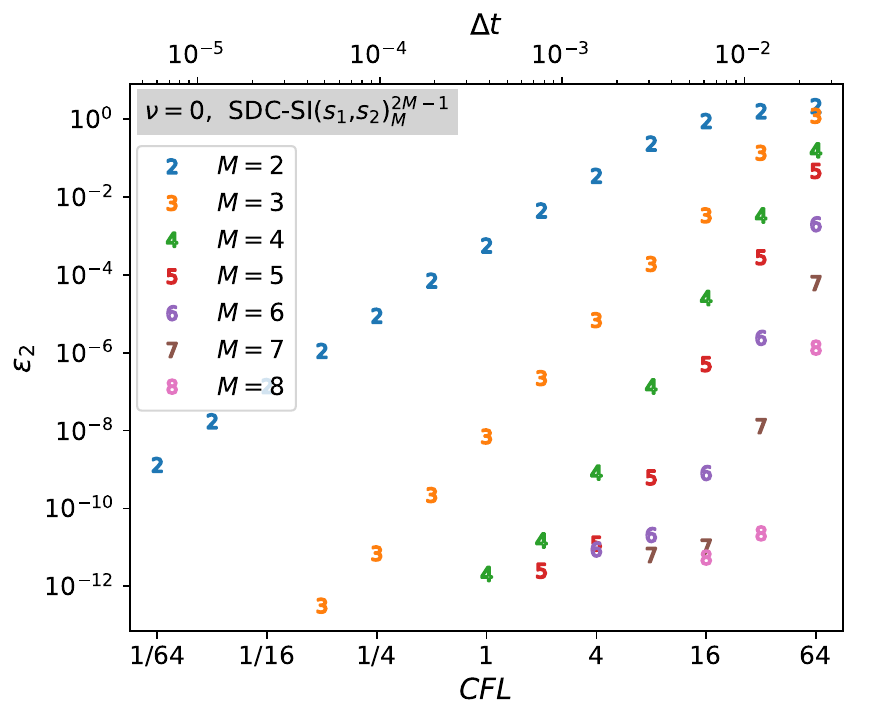}}
  \hfill
  \subcaptionbox{SDC-SI experimental order for ${\nu=0}$
    \label{fig:conv-diff:wave:nu=0e-0:sdc:si1s:order}}
    {\includegraphics[scale=0.52]{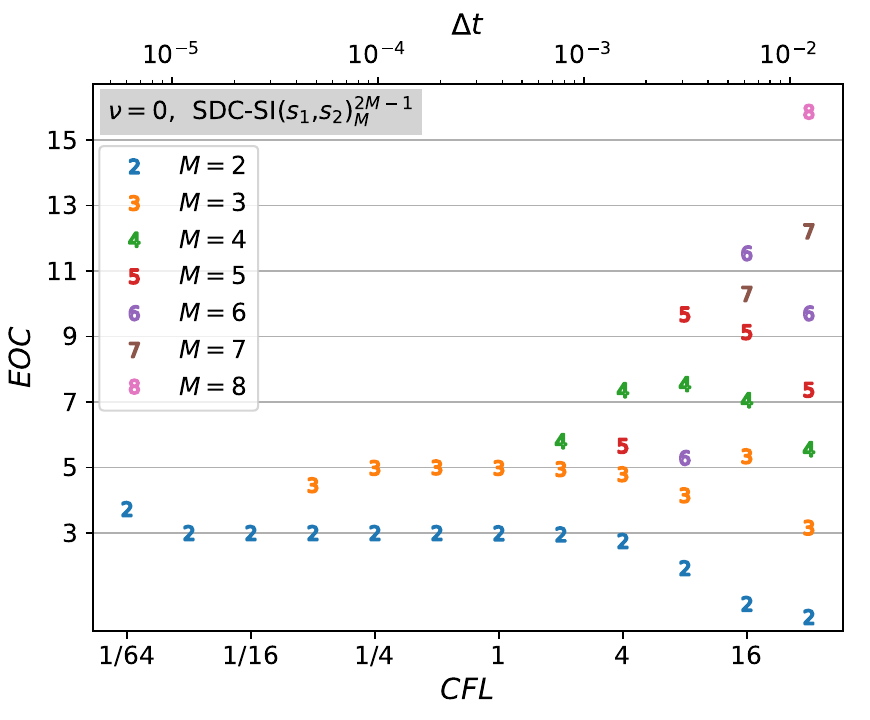}}
  \\[\medskipamount]
  \subcaptionbox{SDC-SI error for ${\nu=10^{-3}}$
    \label{fig:conv-diff:wave:nu=1e-3:sdc:si1s:error}}
    {\includegraphics[scale=0.52]{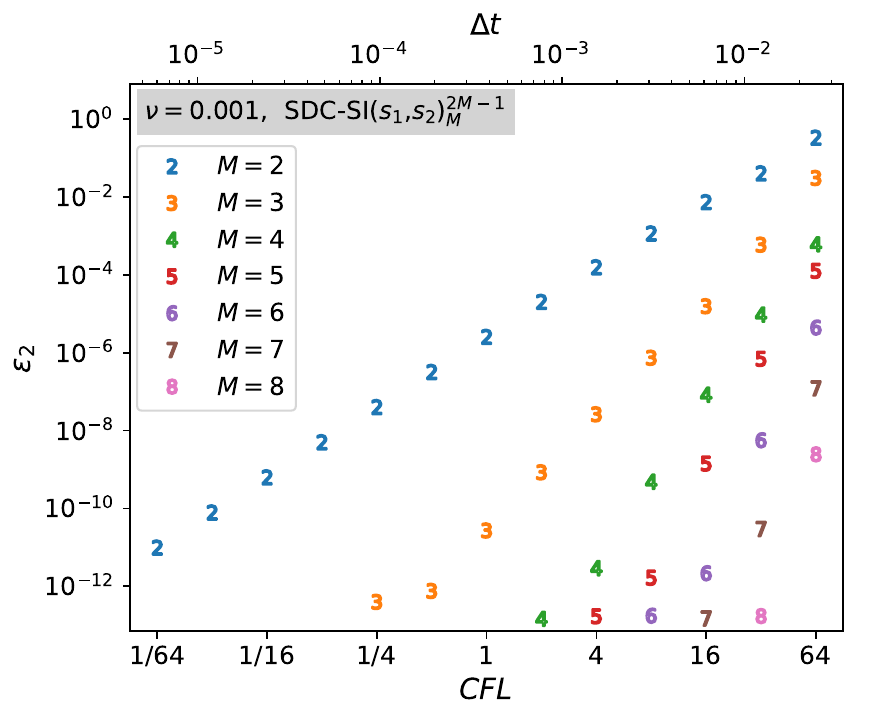}} 
  \hfill
  \subcaptionbox{SDC-SI experimental order for ${\nu=10^{-3}}$
    \label{fig:conv-diff:wave:nu=1e-3:sdc:si1s:order}}
    {\includegraphics[scale=0.52]{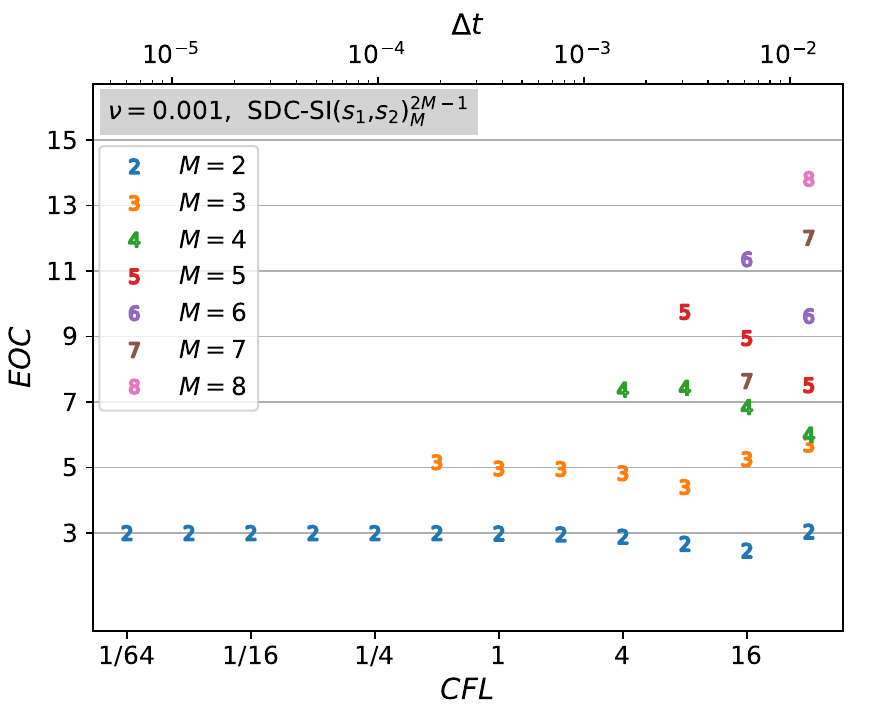}} 
  \caption{Convergence study for the convection-diffusion wave problem.
    \label{fig:conv-diff:wave:overview}}
\end{figure}

The first test considers the case of pure convection in the time interval ${[0,10]}$, i.e. ten periods of the longest wave. 
%%%
Starting from $64$, the CFL number is successively reduced by a factor of 2 down to $1/256$.
%%%
Figure~\ref{fig:conv-diff:wave:nu=0e-0:base:error} shows the $L^2$ error at the final time for selected standalone integrators.
%%%
Due to the absence of diffusivity, SI1(1) and SI2(2) produce virtually identical results.
%%%
They are stable over the entire range and converge at a rate of approximately $\Delta t^2$ for CFL numbers below $2$.
%%%
The second-order IMEX-RK CB2 is less accurate and only succeeds at ${\CFL < 1/16}$ when the numerical diffusion introduced by the Riemann fluxes is large enough to stabilize the method.
%%%
IMEX-RK ARS(4,4,3) and TVD-RK3 achieve their theoretical order of 3 and are 
the most accurate of the standalone integrators tested.
%%%
Note that TVD-RK3 remained stable up to ${\CFL \approx 0.87}$, which is close to the 
theoretical threshold of $1$ and thus confirms the chosen definition of the CFL number.
%%%
Figure~\ref{fig:conv-diff:wave:nu=0e-0:sdc:eu:error} shows the corresponding results for the Euler based SDC methods with ${M=2}$ to 8 Radau points and ${2M-1}$ iterations in every time step.
%%%
With ${M=2}$ or ${M=3}$ the method achieves the expected convergence rate of ${2M-1}$.
%%%
Increasing the number of collocation points improves the accuracy rapidly, while the stability threshold rises slowly from  ${\CFL = 1/2}$ to $2$.
%%%
In contrast, the SDC-SI method with up to 8 Radau points and $s_1$, $s_2$ stages according to Tab.~\ref{tab:sdc:stable:overview} proved stable over the entire $\CFL$ range (Fig.~\ref{fig:conv-diff:wave:nu=0e-0:sdc:si1s:error}).
%%%
With ${M = 2}$ to $4$ points the method achieves the expected order of ${2M-1}$. 
%%%
For larger $M$ the convergence rate cannot be estimated due to the missing asymptotic range (Fig.~\ref{fig:conv-diff:wave:nu=0e-0:sdc:si1s:order}).
%%%

The second test assumes a diffusivity of ${\nu = 0.001}$.
%%%
Figure~\ref{fig:conv-diff:wave:nu=1e-3:sdc:si1s:error} shows a faster convergence in this case, which can be attributed to the damping of the short-wave components and the simpler solution structure at the end time.
%%%
The experimental order of convergence, on the other hand, changes only slightly (Fig.~\ref{fig:conv-diff:wave:nu=1e-3:sdc:si1s:order}).
%%%

\begin{table}
  \begin{tabular}{lccccccc} 
  \toprule
  $i$        & $1$    & $2$    & $3$     & $4$     & $5$     & $6$     & $7$     \\\midrule
  $\kappa_i$ & $2\pi$ & $6\pi$ & $10\pi$ & $14\pi$ & $18\pi$ & $24\pi$ & $30\pi$ \\
  $a_i$      & $1.00$ & $1.50$ & $1.80$  & $1.70$  & $1.50$  & $1.30$  & $1.15$  \\
  $s_i$      & $0.00$ & $0.05$ & $0.10$  & $0.15$  & $0.20$  & $0.30$  & $0.18$  \\
  \bottomrule
  \end{tabular}
  \caption{Wave packet coefficients
    \label{tab:wave-packet}}
\end{table}

%===================================================================================

\subsection{Burgers equation}
\label{sec:numerical-experiments:burgers}

%-----------------------------------------------------------------------------------

The Burgers equation poses a challenge to numerical methods as it allows for nonlinear interactions, including the formation of discontinuities.
%%%
Furthermore, the resulting variation of the convective Jacobian matrix leads to a variable diffusivity in the semi-implicit methods developed in this work.

\begin{figure}
  \subcaptionbox{error
    \label{fig:burgers:wave:sdc:si1s:error}}
    {\includegraphics[scale=0.53]{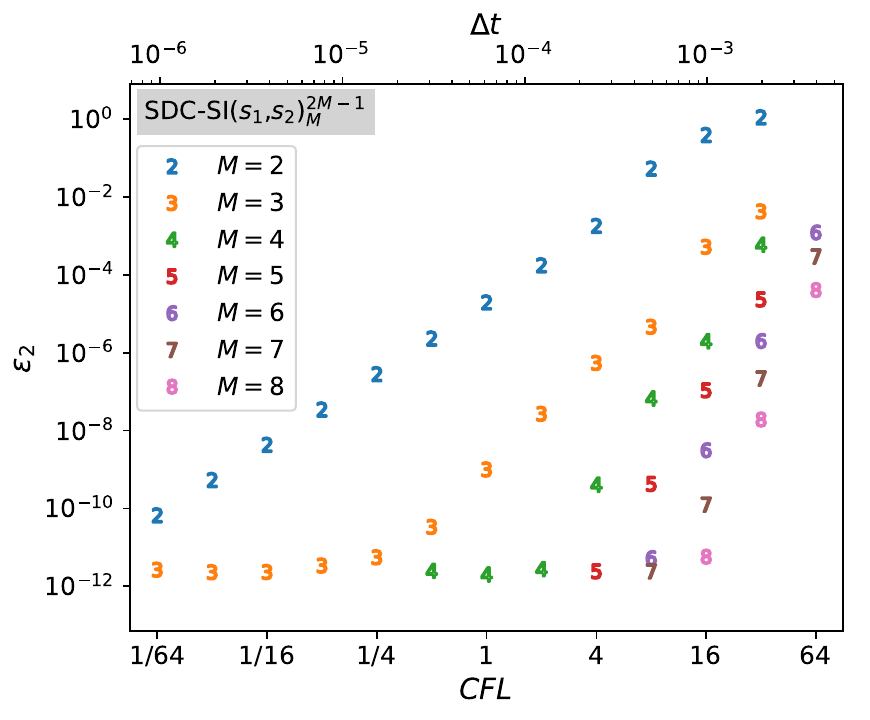}}
  \hfill
  \subcaptionbox{experimental order of convergence
    \label{fig:burgers:wave:sdc:si1s:order}}
    {\includegraphics[scale=0.53]{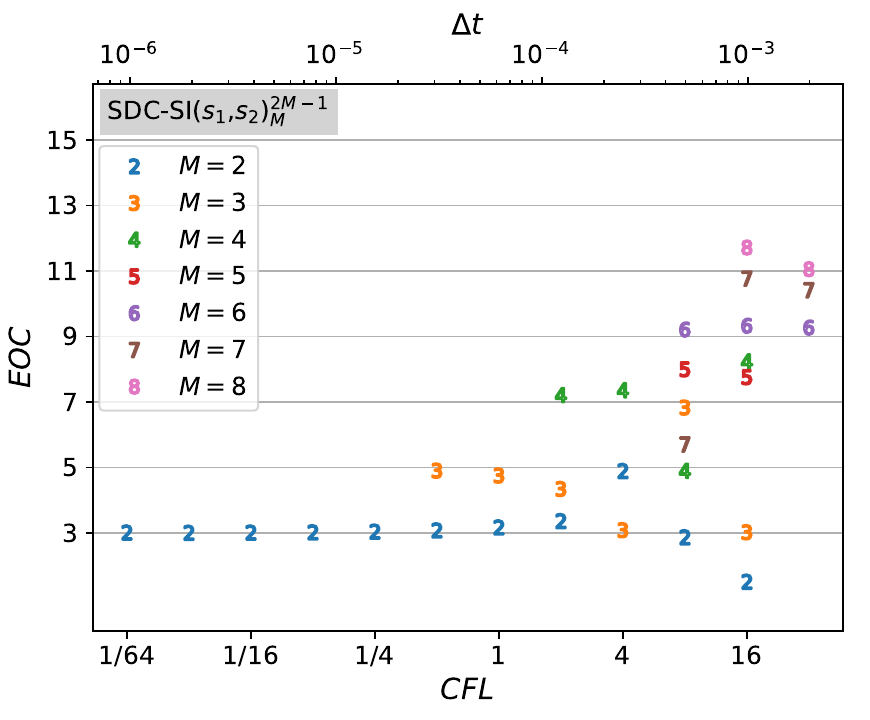}} 
  \caption{%
    Convergence of the SDC-SI method for the inviscid Burgers wave problem 
    \label{fig:burgers:wave}}
\end{figure}

%%%
The first test case examines the wave packet \eqref{eq:wave-packet} with the parameters given in Tab.~\ref{tab:wave-packet}.
%%%
In order to satisfy the Burgers equation, the source term is defined as
\begin{equation*}
  f_\mathrm{s}(x,t) = \d_t u + \d_x \bigl(\tfrac{1}{2}u^2\bigr) - \d_x(\nu \d_x u) 
  \,,
\end{equation*}
based on the exact solution ${u(x,t)}$.
%%%
The computational domain and the discretization parameters are identical to the convection-diffusion problem above.
%%%
As an example, Figure \ref{fig:burgers:wave} shows the results obtained with the SDC-SI method for ${\nu = 0}$.
%%%
Despite the nonlinearity of the Burgers equation, the method proved to be robust and succeeded in all cases with CFL numbers of at least 32.
%%%
The observed errors and convergence rates agree well with the corresponding results for the convection-diffusion problem in Fig.~\ref{fig:conv-diff:wave:nu=0e-0:sdc:si1s:error},\subref{fig:conv-diff:wave:nu=0e-0:sdc:si1s:order}.
%%%

\begin{figure}
  \subcaptionbox{SDC-EU error with marginal spatial resolution
    \label{fig:burgers:front:sdc:eu:error}}
    {\includegraphics[scale=0.53]{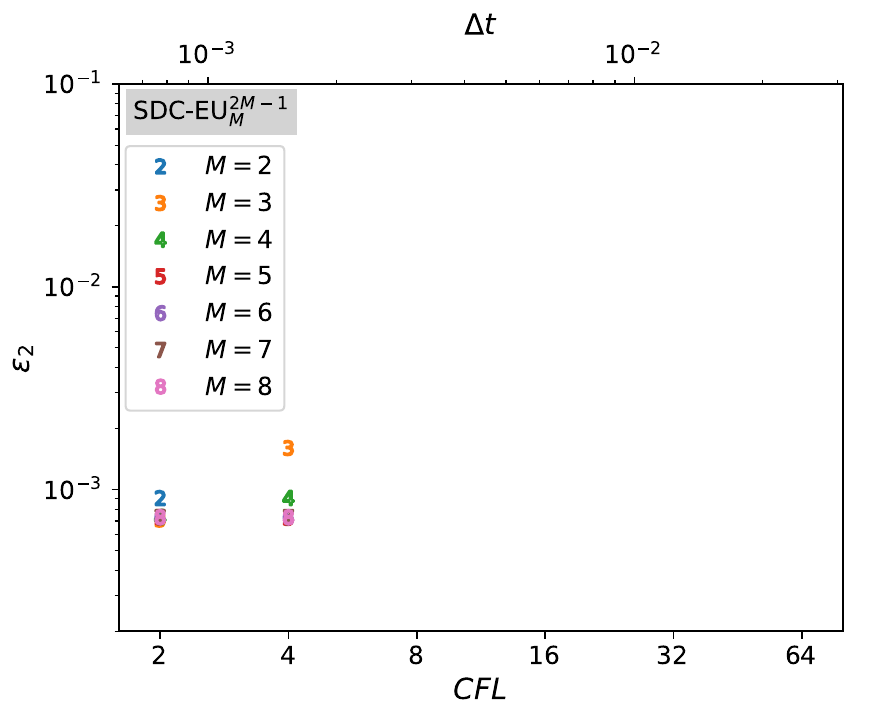}}
  \hfill
  \subcaptionbox{SDC-SI error with marginal spatial resolution
    \label{fig:burgers:front:sdc:si1s:error}}
    {\includegraphics[scale=0.53]{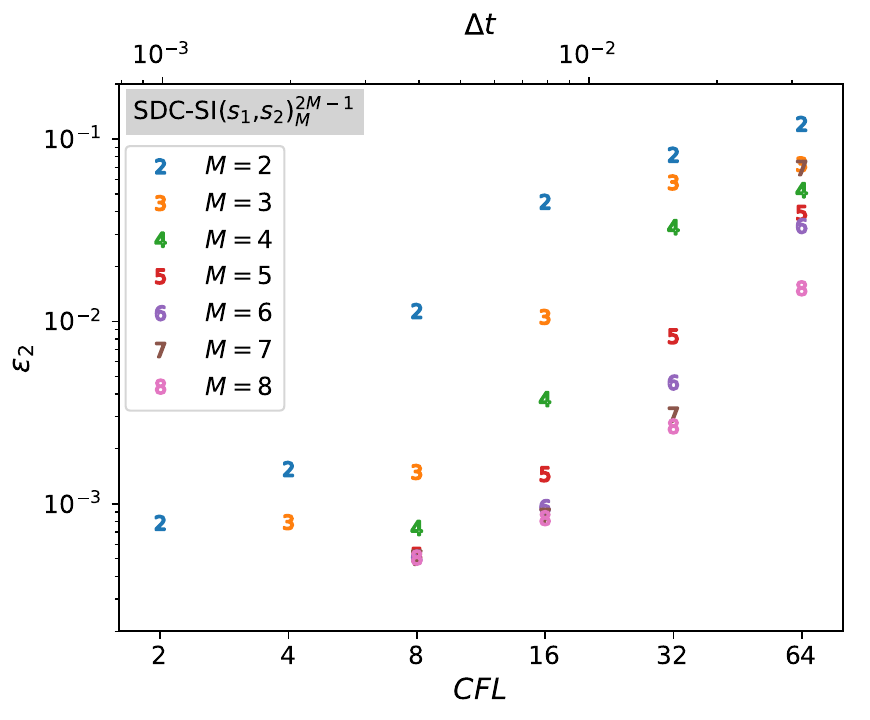}} 
  \\[\medskipamount]
  \subcaptionbox{SDC-SI solution with marginal spatial resolution
    \label{fig:burgers:front:sdc:si1s:cfl:hires}}
    {\includegraphics[scale=0.53]{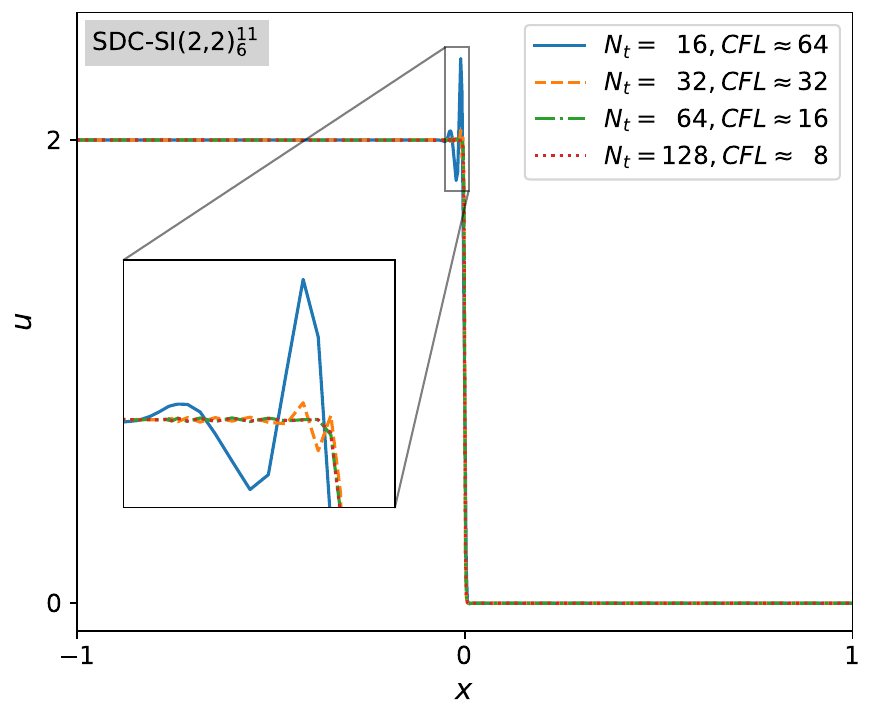}}
  \hfill
  \subcaptionbox{SDC-SI solution with low spatial resolution and \linebreak 
    artificial diffusion
    \label{fig:burgers:front:sdc:si1s:cfl:lores}}
    {\includegraphics[scale=0.53]{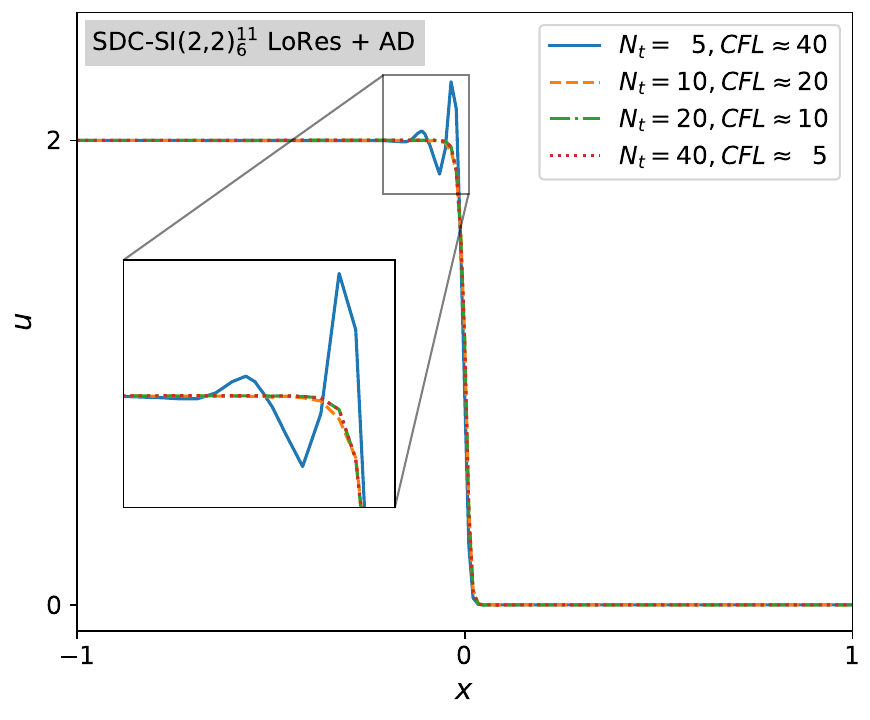}}
  \caption{%
    Results for the Burgers moving front problem
    \label{fig:burgers:front}}
\end{figure}

The next test case relates to the moving front solution
\begin{equation*}
  u(x,t) = 1 - \tanh\left(\frac{x + 0.5 - t}{2 \nu}\right)
\end{equation*}
for ${x \in [-1,1]}$, see \cite{SE_Hesthaven2008a}.
%%%
Assuming a diffusivity of $\nu = 0.001$ yields a thin front advancing with unit velocity.
%%%
Two scenarios are considered: 
the first one uses 50 elements of degree 15 and 
the second 20 elements of degree 10.
%%%
Time integration is performed over the interval ${[0,0.5]}$ with $N_t$ steps.
%%%
The first scenario is characterized by a marginal spatial resolution, which is sufficient to capture the front with an $L^2$ error of about ${6 \times 10^{-4}}$.
%%%
Figure~\ref{fig:burgers:front:sdc:eu:error} shows the corresponding convergence results for the Euler based SDC-EU method as a basis for comparison.
%%%
Note that the method achieves optimal accuracy with $\CFL \le 4$, but fails with larger steps.
%%%
In contrast, the SDC-SI method succeeds for all $M$ with only ${N_t=16}$ steps or ${\CFL \approx 64}$ (Fig.~\ref{fig:burgers:front:sdc:si1s:error}).
%%%
As the CFL number decreases, the method converges and reaches the minimal error with half as many steps as SDC-EU.
%%%
Figure~\ref{fig:burgers:front:sdc:si1s:cfl:hires} shows the behavior of the solution obtained with SDC-SI($2$,$2$)$_6^{11}$ for different CFL numbers.
%%%
Although spurious oscillations occur, the method remains stable up to ${\CFL = 64}$.
%%%
The inset gives a close view on the upstream region near the front and illustrates how the oscillations reduce as the time step decreases.
%%%
At ${\CFL \le 16}$ the minimum error is reached and leaves slight ripples due to the marginal spatial resolution.
%%%

In the second scenario the artificial diffusivity \eqref{eq:art-diffusivity} is used with parameters ${\ds = 2}$ and ${\cs = 0.4}$.
%%%
Figure~\ref{fig:burgers:front:sdc:si1s:cfl:lores} shows the numerical solution obtained with SDC-SI($2$,$2$)$_6^{11}$ for selected step numbers $N_t$.
%%%
The method succeeds already with 4 steps, although with oscillations somewhat larger as depicted for ${N_t = 5}$.
%%%
As the number of steps increases, these oscillations disappear gradually and nearly vanish for ${N_t = 10}$ or ${\CFL \approx 20}$.
%%%
The solutions obtained with 20 and 40 steps are close to the temporally converged case and hence almost identical.
%%%

%===================================================================================

\subsection{Euler and Navier-Stokes equations}
\label{sec:numerical-experiments:euler+cns}

%-----------------------------------------------------------------------------------

The final series of tests deals with a selection of compressible flow problems.
%%%
The fluid is a perfect gas with ${R = 287.28}$ and ${\gamma = 1.4}$, which corresponds to air when assuming SI units.
%%%

\begin{figure}
  \subcaptionbox{${M = 3}$ Radau points
    \label{fig:cns:wave:rr3}}
    {\includegraphics[scale=0.52]{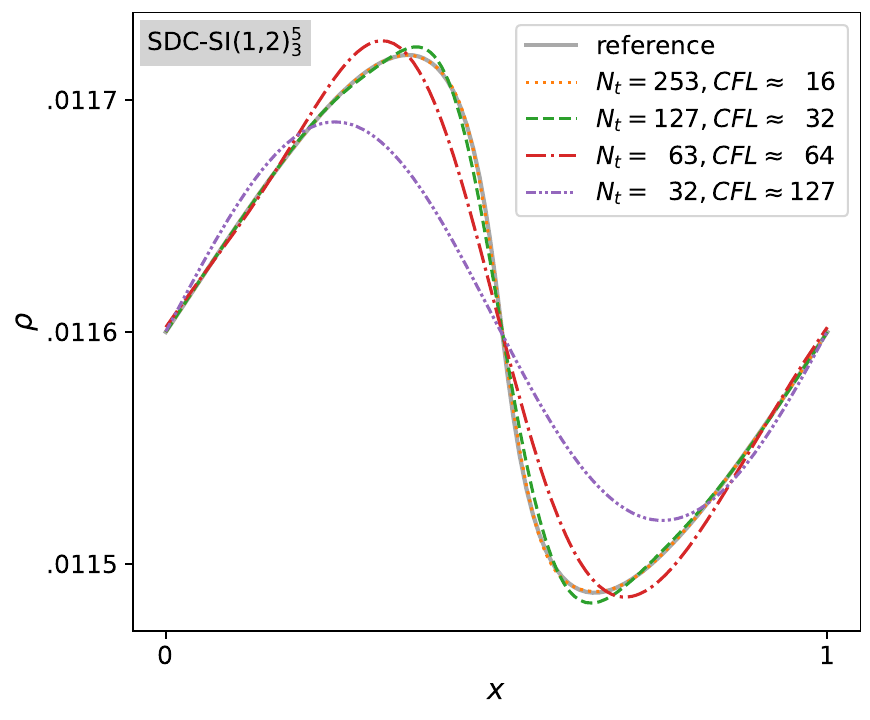}}
  \hfill
  \subcaptionbox{${M = 8}$ Radau points
    \label{fig:cns:wave:rr8}}
    {\includegraphics[scale=0.52]{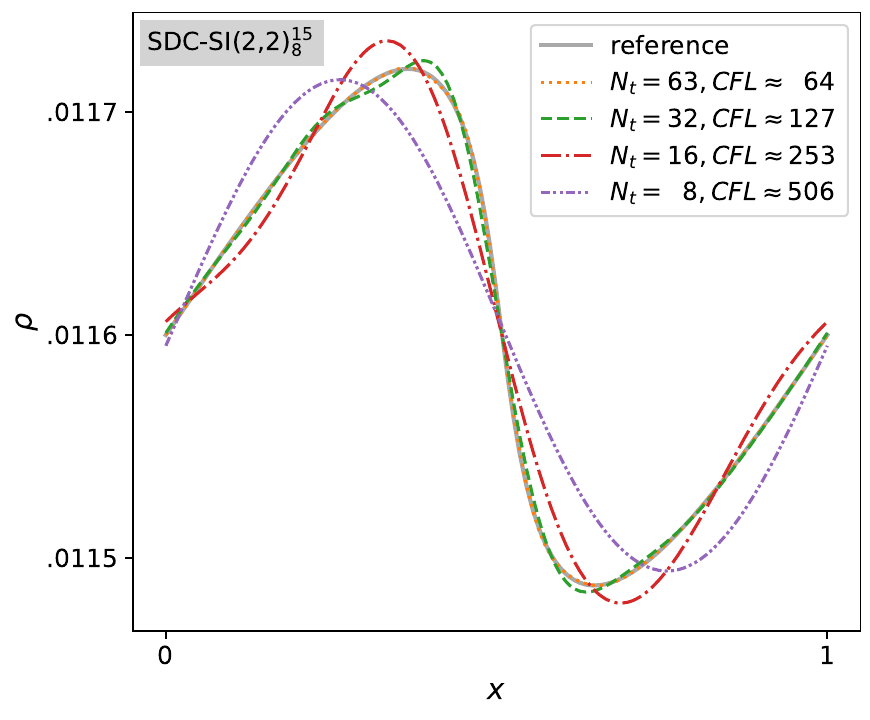}} 
  \caption{%
    Acoustic wave density distributions 
    obtained with SDC-SI($s_1$,$s_2$)$_{M}^{2M-1}$
    \label{fig:cns:wave}}
\end{figure}

The first problem describes an acoustic wave traveling downstream in a uniform flow with 
Mach number ${\mathit{Ma} = 0.1}$, 
pressure ${p_{\infty} = 1000}$,
temperature ${T_{\infty} = 300}$,
dynamic viscosity ${\eta = 10^{-5}}$ and
Prandtl number ${\mathit{Pr} = 0.75}$.
%%%
It is initialized by superimposing the periodic velocity perturbation
\begin{equation*}
  v'(x) = \hat v \sin(2 \pi x)
\end{equation*}
with an amplitude of 
${\hat v = 10^{-2} a_{\infty}}$, where 
${a_{\infty} \approx 347.36}$ is the speed of sound in the undisturbed fluid.
%%%
Assuming an isotropic state, the initial temperature is given by
\begin{equation*}
  T_0(x) = \frac{1}{\gamma R}\left( a_{\infty} +\frac{\gamma-1}{2}v' \right)^2
  \,,
\end{equation*}
the pressure by
\begin{equation*}
  p_0(x) = p_{\infty} \left(\frac{T}{T_{\infty}}\right)^{\gamma/(\gamma-1)}
\end{equation*}
and the velocity by
\begin{equation*}
  v_0(x) = \mathit{Ma} \, a_{\infty} + v'(x)
  \,.
\end{equation*}
%%%
Note that the amplitude is large enough to expect a nonlinear steepening.
%%%
In contrast, viscous damping is negligible due to the high Reynolds number 
${\mathit{Re} = \rho_{\infty} a_{\infty} / \eta \approx 4 \times 10^5}$.
%%%
The evolution of the wave is studied in the periodic domain ${\Omega = [0,1]}$ from ${t = 0}$ to ${t_{\mathrm{end}} = 2.619 \times 10^{-2}}$, which corresponds to approximately 10 convective units.
%%%
For spatial discretization, the domain is decomposed into 10 elements of degree 15.
%%%
Since no exact solution is available, the TVD-RK3 method with ${\Delta t \approx 10^{-7}}$ or ${\CFL \approx 0.016}$ serves as a reference.
%%%
The SDC-SI method was applied with ${M=2}$ to 8 Radau points and ${2M-1}$ iterations.
%%%
Between ${N_t = 4}$ and $64,826$ time steps were used for time integration, which corresponds to CFL numbers in the range from $1/32$ to $1013$.
%%%
The fact that all runs were successful underlines the stability of the method.
%%%
As an example, the computed density distributions are shown in
Fig.~\ref{fig:cns:wave:rr3} for ${M = 3}$ and in
Fig.~\ref{fig:cns:wave:rr8} for ${M = 8}$.
%%%
Note that the results are visually indistinguishable from the reference for 
${\CFL \approx 16}$ in the first case and 
${\CFL \approx 64}$ in the second case.
%%%
Hence, only about $6$ time steps with $8$ subintervals are required for one convective unit.
%%%

\begin{figure}
  \subcaptionbox{density
    \label{fig:cns:shock-tube:rho}}
    {\includegraphics[scale=0.52]{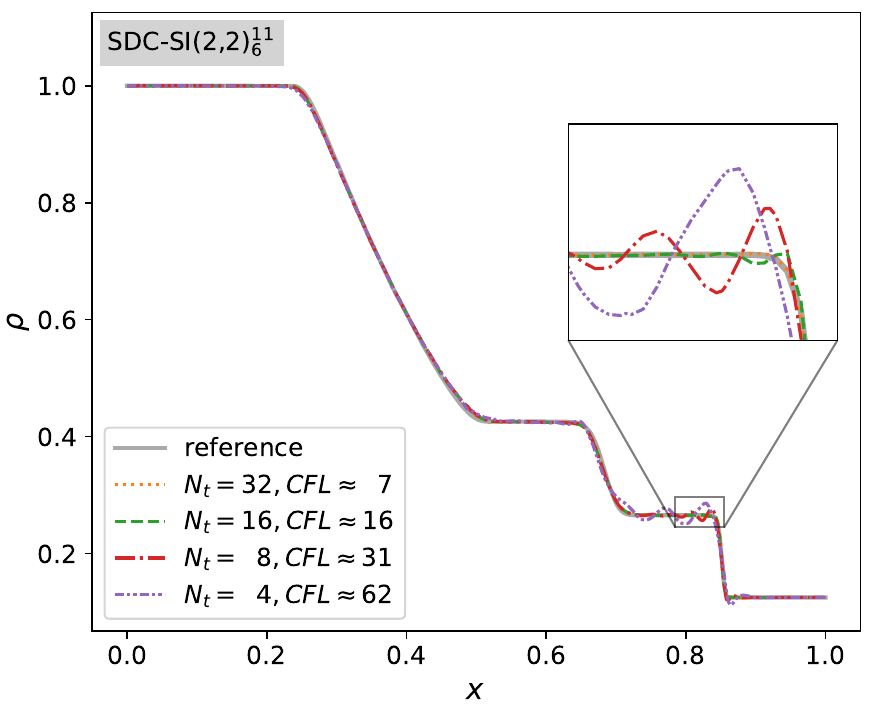}}
  \hfill
  \subcaptionbox{velocity
    \label{fig:cns:shock-tube:v}}
    {\includegraphics[scale=0.52]{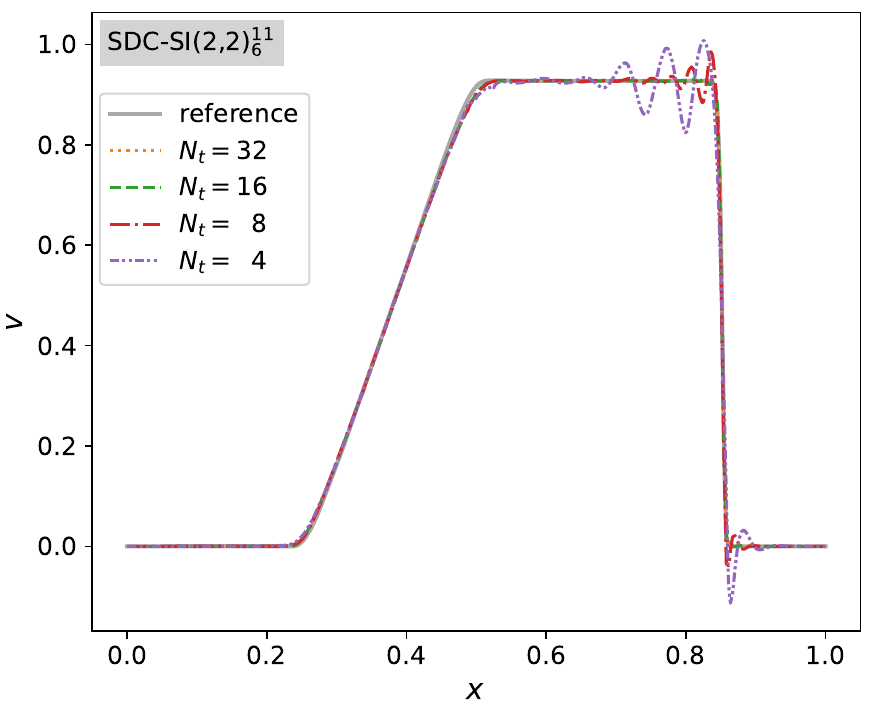}} 
  \caption{%
    Shock tube results obtained with 
    SDC-SI SDC-SI($2$,$2$)$_{6}^{11}$ 
    for different CFL numbers
    \label{fig:cns:shock-tube}}
\end{figure}

The second example is Sod's shock tube problem for the Euler equations with initial conditions \cite{BM_Sod1978a}
%\begin{subequations}
  \begin{alignat*}{4}
    &\rho = 1,     &\quad& p = 1,   &\quad& v = 0 && \qquad\text{if} \quad x <   0.5 \\
    &\rho = 0.125, &\quad& p = 0.1, &\quad& v = 0 && \qquad\text{if} \quad x \ge 0.5
  \end{alignat*}
%\end{subequations}
to solve on ${\Omega = [0,1]}$ until ${t = 0.2}$.
%%%
For the numerical studies, the domain is divided into 80 elements of degree 5.
%%%
The artificial diffusivity \eqref{eq:art-diffusivity} is used with ${\ds = 6}$ and ${\cs = 0.4}$ to capture the evolving shock and the contact discontinuity.
%%%
Figure~\ref{fig:cns:shock-tube} shows the density and velocity distributions obtained with SDC-SI($2$,$2$)$_6^{11}$ at different step sizes.
%%%
TVD-RK3 with $10^5$ steps or ${\Delta t = 2\times 10^{-6}}$ serves as a reference.
%%%
Note that the SDC-SI method produces reasonable results with only 4 steps.
%%%
This is remarkable, because the shock wave traverses a quarter of the domain during this time.
%%%
Although the density and velocity distributions behind the shock show considerable oscillations, these decrease as the number of steps increases.
%%%
With ${N_t = 32}$ or ${\CFL \approx 7}$ they disappear and the solution is visually identical to the reference.
%%%

\begin{figure}
  \subcaptionbox{variation of step size
    \label{fig:cns:shu-osher:cfl}}
    {\includegraphics[scale=0.52]{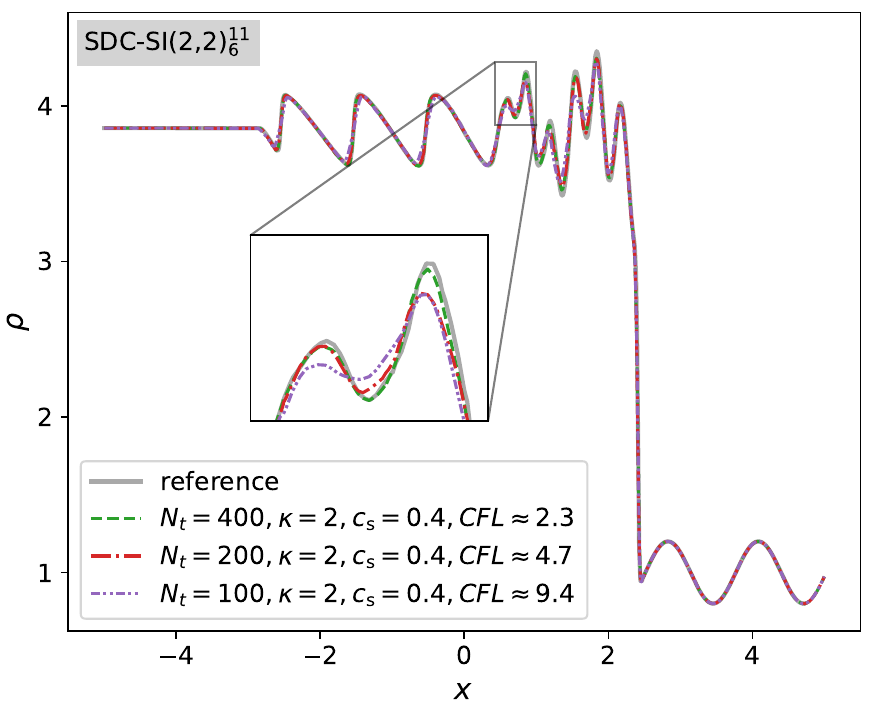}}
  \hfill
  \subcaptionbox{variation of artificial diffusion parameters
    \label{fig:cns:shu-osher:ad}}
    {\includegraphics[scale=0.52]{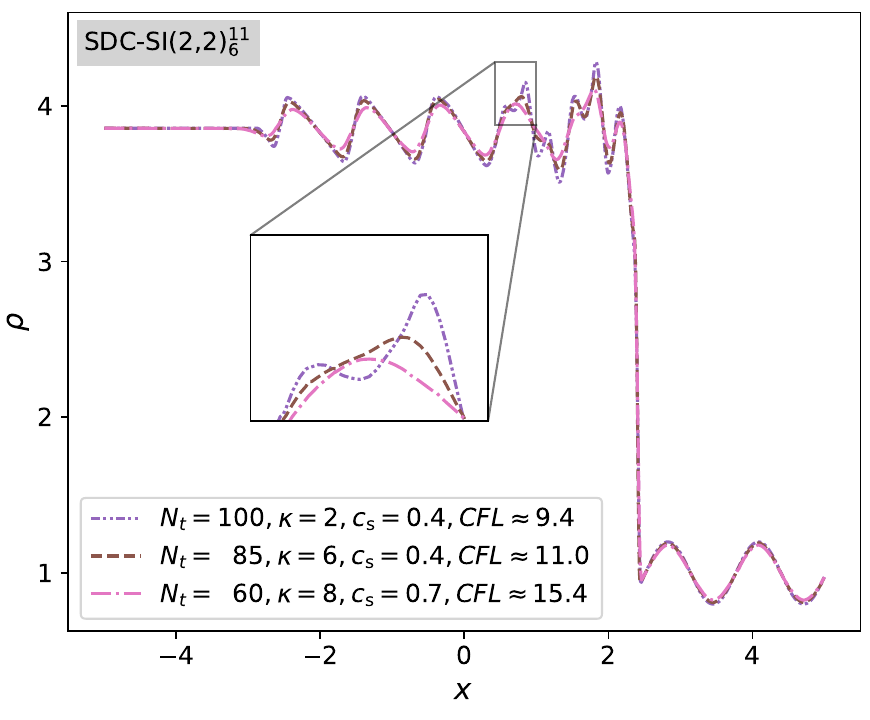}}
  \caption{%
    Shock-fluctuation interaction: density distributions obtained with
    SDC-SI($2$,$2$)$_{6}^{11}$ for different step sizes (left) and 
    different artificial diffusion parameters (right)
    \label{fig:cns:shu-osher}}
\end{figure}

The final test case is the shock-fluctuation benchmark of \citet[Example 8]{TI_Shu1989a}
for the Euler equations with initial conditions
%\begin{subequations}
  \begin{alignat*}{4}
    &  \rho = 3.857143, &\quad& p = 10.33333, &\quad& v = 2.629369 
    && \qquad \text{if} \quad x < -4 \,,\\
    &  \rho = 1 + 0.2\sin 5x, &\quad& p = 1,  &\quad& v = 0 
    && \qquad \text{if} \quad x \ge -4 \,,
  \end{alignat*}
%\end{subequations}
to solve on ${\Omega = [-5, 5]}$ until ${t = 1.8}$.
%%%
This case is more challenging than Sod's problem as it involves the complex interaction between a strong discontinuity and smooth wave-like pattern.
%%%
The numerical solutions were computed with 160 elements of degree 5 using the artificial diffusivity \eqref{eq:art-diffusivity} for stabilization.
%%%
Figure~\ref{fig:cns:shu-osher:cfl} shows the density distributions at ${t = 1.8}$ for SDC-SI($2$,$2$)$_{6}^{11}$ with varying temporal resolution and fixed stabilization parameters ${\ds = 2}$, ${\cs = 0.4}$.
%%%
The given CFL numbers were determined after ${N_t / 2}$ steps.
%%%
For ${N_t \ge 400}$ or ${\CFL \lesssim 2.3}$ the results agree well with reference data obtained with TVD-RK3 using $180,000$ steps or ${\Delta t = 1\times 10^{-5}}$ .
%%%
With increasing step size, the error in the interaction zone behind the shock grows,  but remains at a moderate level up to ${N_t = 100}$ or ${\CFL \approx 9.3}$.
%%%
Using less than 90 steps leads to negative density or temperature values in the intermediate solution and thus to failure of the method.
%%%
In further tests, the SDC-SI method achieved CFL numbers between 8 with ${M=2}$ Radau points and 21 with ${M=16}$.
%%%
A closer examination revealed that the resolution near the stability threshold is insufficient to capture the temporal development.
%%%
As a consequence, the temporal interpolation leads to spurious oscillations which result in physically invalid intermediate values.
%%%
In agreement with this observation, the stability threshold grows if the solution is smoothed by increasing the artificial diffusion.
%%%
Figure~\ref{fig:cns:shu-osher:ad} gives an example for SDC-SI($2$,$2$)$_{6}^{11}$, where the use of a wider sensor interval and a higher maximum diffusivity increased the admissible step size by a factor of about $1.7$.
%%%

%===================================================================================

\subsection{Computational efficiency}
\label{sec:numerical-experiments:efficiency}

%-----------------------------------------------------------------------------------

\begin{figure}
  \subcaptionbox{pure convection
    \label{fig:conv-diff:wave:nu=0e-0:error:t_run}}
    {\includegraphics[scale=0.52]{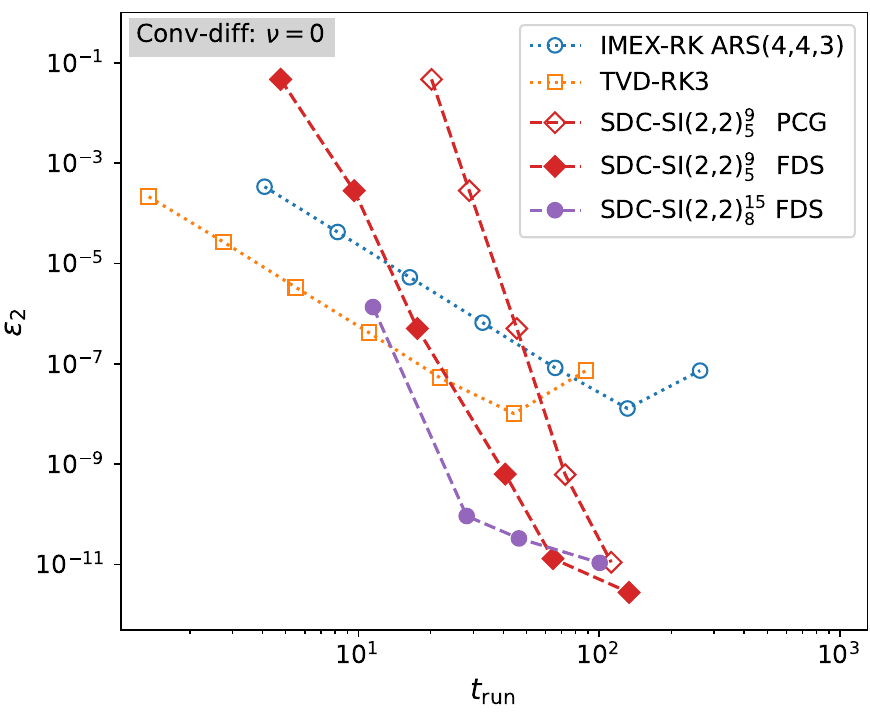}}
  \hfill
  \subcaptionbox{convection-diffusion
    \label{fig:conv-diff:wave:nu=1e-1:error:t_run}}
    {\includegraphics[scale=0.52]{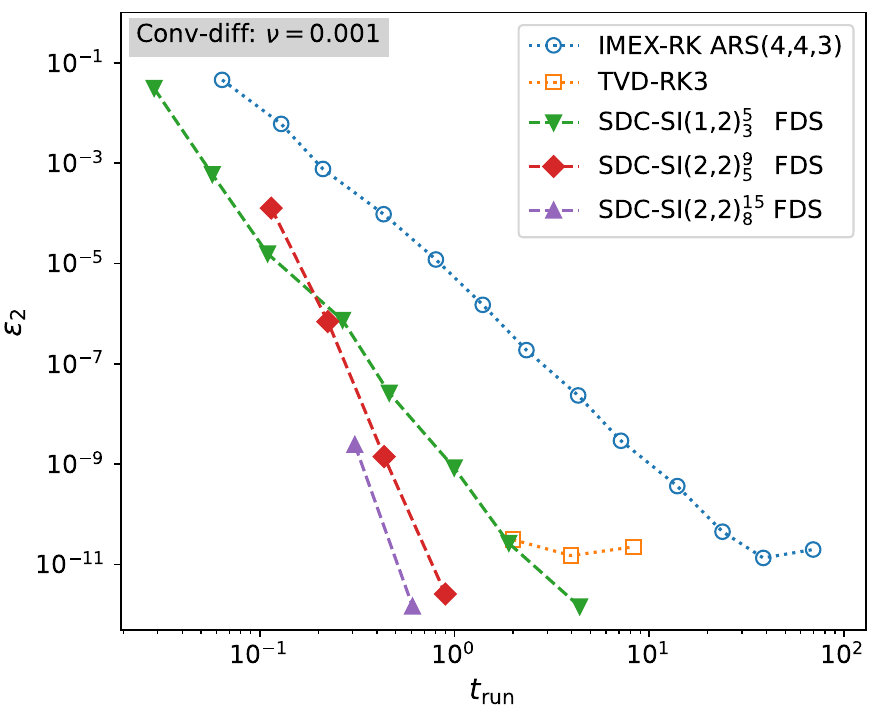}}
  \caption{%
    Error-runtime comparison between Runge-Kutta and SDC methods using either 
    the preconditioned conjugate gradient method (PCG) or the fast direct solver (FDS)
    \label{fig:conv-diff:wave:error:t_run}}
\end{figure}

In the numerical experiments described above, the focus was deliberately placed on stability and accuracy.
%%%
The iteration counts and terminal conditions for solving the implicit equation systems were chosen carefully to exclude adverse effects due to algebraic errors.
%%%
In addition, the iterative solvers have not yet been optimized and are therefore inferior to explicit methods such as TVD-RK3.
%%%
However, to demonstrate the potential of the proposed SDC methods, performance tests were performed for the convection-diffusion equation for which a fast direct solver is available.
%%%
Figure~\ref{fig:conv-diff:wave:nu=0e-0:error:t_run} shows a comparison between  Runge-Kutta and SDC methods for the purely convective wave problem studied in Sec.~\ref{sec:numerical-experiments:conv-diff}.
%%%
For errors down to $10^{-6}$ TVD-RK3 outperforms the other methods, which is not surprising because the method was designed for this type of problems.
%%%
IMEX-RK ARS(4,4,3) requires about three times the time to achieve the same error.
%%%
The SDC methods are more expensive in this error range, but become competitive with increasing accuracy, especially when using the fast solver.
%%%
For example, the SDC-SI(2,2) method with 8 Radau points surpasses TVD-RK3 for errors below $10^{-6}$.
%%%
In presence of diffusion, the advantage of explicit Runge-Kutta methods diminishes further.
%%%
As an example, Figure~\ref{fig:conv-diff:wave:nu=1e-1:error:t_run} shows the results for ${\nu = 0.001}$.
%%%
Although convection still dominates, the SDC methods with 3, 5 or 8 Radau points clearly outperform the Runge-Kutta methods in the entire error range.

%===================================================================================

% !TEX encoding = UTF-8 Unicode
% !TEX root = ssi-sdc-1d.tex

%%%%%%%%%%%%%%%%%%%%%%%%%%%%%%%%%%%%%%%%%%%%%%%%%%%%%%%%%%%%%%%%%%%%%%%%%%%%%%%%%%%%

\section{Conclusions}
\label{sec:conclusions}

%===================================================================================

This work presents a new family of spectral deferred correction (SDC) methods with excellent stability properties.
The key ingredient are novel time integrators which combine a semi-implicit splitting of the convection term inspired by the Lax-Wendroff method with an implicit treatment of diffusion and sources.
%%%
These integrators possess excellent stability properties and require only the solution of positive definite or semi-definite linear systems.
%%%
In the present study they were used to construct semi-implicit SDC methods on the basis of right-sided Radau points.
%%%
Numerical evidence suggests that these methods are $L$-stable with 2 to 6 points or orders from 3 to 11 and require very little diffusion for orders 13 and 15.
%%%

In numerical experiments, the SDC methods were applied to high-order discontinuous Galerkin formulations of one-dimensional conservation laws.
%%%
Tests with the convection-diffusion equation proved the excellent stability and accuracy of the methods: They remained stable with CFL numbers up to at least 64 and achieved the expected order of convergence, i.e. ${2M-1}$ with $M$ Radau points.
%%%
In comparison, the SDC method based on the IMEX Euler scheme reached CFL numbers of only 0.5 with ${M=2}$ and 2 with ${M \ge 4}$.
%%%
Excellent results were also obtained for nonlinear problems, including the Burgers, Euler and Navier-Stokes equations.
%%%
In particular, it could be shown that the proposed SDC methods can handle marginally resolved solutions and easily accommodate artificial diffusion techniques for stabilization near discontinuities.
%%%
Future work should focus on acceleration, e.g., by using multilevel techniques and generalization to multidimensional problems \cite{TI_Speck2015a}.
%%%

%===================================================================================

\section*{Data Availability}

The datasets generated during the current study are available from the corresponding author on reasonable request.

\section*{Funding}

Funding by German Research Foundation (DFG) in frame of the project STI 57/9-1 is gratefully acknowledged. The author would like to thank ZIH, TU Dresden, for the provided computational resources. Open Access funding is enabled and organized by Projekt DEAL.

\section*{Ethics declarations}

The author declares that he has no known competing financial interests or personal relationships that could have appeared to influence the work reported in this paper.

%%%%%%%%%%%%%%%%%%%%%%%%%%%%%%%%%%%%%%%%%%%%%%%%%%%%%%%%%%%%%%%%%%%%%%%%%%%%%%%%%%%%
% Bibliography


\begin{thebibliography}{67}
\providecommand{\natexlab}[1]{#1}
\providecommand{\url}[1]{\texttt{#1}}
\expandafter\ifx\csname urlstyle\endcsname\relax
  \providecommand{\doi}[1]{doi: #1}\else
  \providecommand{\doi}{doi: \begingroup \urlstyle{rm}\Url}\fi

\bibitem[Ascher et~al.(1995)Ascher, Ruuth, and Wetton]{TI_Ascher1995a}
U.~M. Ascher, S.~J. Ruuth, and B.~T.~R. Wetton.
\newblock Implicit-explicit methods for time-dependent partial differential
  equations.
\newblock \emph{SIAM Journal on Numerical Analysis}, 32\penalty0 (3):\penalty0
  797--823, 1995.
\newblock \doi{10.1137/0732037}.

\bibitem[Ascher et~al.(1997)Ascher, Ruuth, and Spiteri]{TI_Ascher1997a}
U.~M. Ascher, S.~J. Ruuth, and R.~J. Spiteri.
\newblock Implicit-explicit {R}unge-{K}utta methods for time-dependent partial
  differential equations.
\newblock \emph{Applied Numerical Mathematics}, 25\penalty0 (2-3):\penalty0
  151--167, 1997.
\newblock \doi{10.1016/s0168-9274(97)00056-1}.

\bibitem[Baeza et~al.(2020)Baeza, B{\"u}rger, Mart{\'\i}, Mulet, and
  Zor{\'\i}o]{TI_Baeza2020a}
A.~Baeza, R.~B{\"u}rger, M.~d.~C. Mart{\'\i}, P.~Mulet, and D.~Zor{\'\i}o.
\newblock On approximate implicit taylor methods for ordinary differential
  equations.
\newblock \emph{Computational and Applied Mathematics}, 39\penalty0 (4), 2020.
\newblock \doi{10.1007/s40314-020-01356-8}.

\bibitem[Barter and Darmofal(2010)]{SE_Barter2010a}
G.~E. Barter and D.~L. Darmofal.
\newblock Shock capturing with pde-based artificial viscosity for dgfem: Part
  i. formulation.
\newblock \emph{Journal of Computational Physics}, 229\penalty0 (5):\penalty0
  1810--1827, 2010.
\newblock \doi{10.1016/j.jcp.2009.11.010}.

\bibitem[Bassi et~al.(2015)Bassi, Botti, Colombo, Ghidoni, and
  Massa]{TI_Bassi2015a}
F.~Bassi, L.~Botti, A.~Colombo, A.~Ghidoni, and F.~Massa.
\newblock Linearly implicit {R}osenbrock-type {R}unge--{K}utta schemes applied
  to the discontinuous {G}alerkin solution of compressible and incompressible
  unsteady flows.
\newblock \emph{Computers {\&} Fluids}, 118:\penalty0 305--320, 2015.
\newblock \doi{10.1016/j.compfluid.2015.06.007}.

\bibitem[Boscarino et~al.(2016)Boscarino, Filbet, and Russo]{TI_Boscarino2016a}
S.~Boscarino, F.~Filbet, and G.~Russo.
\newblock High order semi-implicit schemes for time dependent partial
  differential equations.
\newblock \emph{Journal of Scientific Computing}, 68\penalty0 (3):\penalty0
  975--1001, 2016.
\newblock \doi{10.1007/s10915-016-0168-y}.

\bibitem[Boscheri and Tavelli(2022)]{TI_Boscheri2022a}
W.~Boscheri and M.~Tavelli.
\newblock High order semi-implicit schemes for viscous compressible flows in
  3d.
\newblock \emph{Applied Mathematics and Computation}, 434:\penalty0 127457,
  2022.
\newblock \doi{10.1016/j.amc.2022.127457}.

\bibitem[Brooks and Hughes(1982)]{FE_Brooks1982a}
A.~N. Brooks and T.~J. Hughes.
\newblock Streamline upwind/petrov-galerkin formulations for convection
  dominated flows with particular emphasis on the incompressible navier-stokes
  equations.
\newblock \emph{Computer Methods in Applied Mechanics and Engineering},
  32\penalty0 (1--3):\penalty0 199--259, Sept. 1982.
\newblock \doi{10.1016/0045-7825(82)90071-8}.

\bibitem[Canuto et~al.(2007)Canuto, Hussaini, Quarteroni, and
  Zang]{SE_Canuto2007a}
C.~Canuto, M.~Y. Hussaini, A.~Quarteroni, and T.~A. Zang.
\newblock \emph{Spectral Methods. Evolution to Complex Geometries and
  Applications to Fluid Dynamics}.
\newblock Springer-Verlag GmbH, 2007.
\newblock ISBN 3540307273.

\bibitem[Canuto et~al.(2011)Canuto, Hussaini, Quarteroni, and
  Zang]{SE_Canuto2011a}
C.~Canuto, M.~Y. Hussaini, A.~Quarteroni, and T.~A. Zang.
\newblock \emph{Spectral Methods. Fundamentals in Single Domains}.
\newblock Springer Berlin Heidelberg, 2011.
\newblock ISBN 3540307257.

\bibitem[Causley and Seal(2019)]{TI_Causley2019a}
M.~Causley and D.~Seal.
\newblock On the convergence of spectral deferred correction methods.
\newblock \emph{Communications in Applied Mathematics and Computational
  Science}, 14\penalty0 (1):\penalty0 33--64, 2019.
\newblock \doi{10.2140/camcos.2019.14.33}.

\bibitem[Cavaglieri and Bewley(2015)]{TI_Cavaglieri2015a}
D.~Cavaglieri and T.~Bewley.
\newblock Low-storage implicit/explicit {R}unge-{K}utta schemes for the
  simulation of stiff high-dimensional {ODE} systems.
\newblock \emph{Journal of Computational Physics}, 286:\penalty0 172--193,
  2015.
\newblock \doi{10.1016/j.jcp.2015.01.031}.

\bibitem[Christlieb et~al.(2009)Christlieb, Ong, and Qiu]{TI_Christlieb2009a}
A.~Christlieb, B.~Ong, and J.-M. Qiu.
\newblock Comments on high-order integrators embedded within integral deferred
  correction methods.
\newblock \emph{Communications in Applied Mathematics and Computational
  Science}, 4\penalty0 (1):\penalty0 27--56, 2009.
\newblock \doi{10.2140/camcos.2009.4.27}.

\bibitem[Christlieb et~al.(2011)Christlieb, Morton, Ong, and
  Qiu]{TI_Christlieb2011a}
A.~Christlieb, M.~Morton, B.~Ong, and J.-M. Qiu.
\newblock Semi-implicit integral deferred correction constructed with additive
  {R}unge--{K}utta methods.
\newblock \emph{Communications in Mathematical Sciences}, 9\penalty0
  (3):\penalty0 879--902, 2011.
\newblock \doi{10.4310/cms.2011.v9.n3.a10}.

\bibitem[Christlieb et~al.(2015)Christlieb, Liu, and Xu]{TI_Christlieb2015a}
A.~J. Christlieb, Y.~Liu, and Z.~Xu.
\newblock High order operator splitting methods based on an integral deferred
  correction framework.
\newblock \emph{Journal of Computational Physics}, 294:\penalty0 224--242,
  2015.
\newblock \doi{10.1016/j.jcp.2015.03.032}.

\bibitem[Deuflhard and Bornemann(2002)]{TI_Deuflhard2002a}
P.~Deuflhard and F.~Bornemann.
\newblock \emph{Scientific Computing with Ordinary Differential Equations}.
\newblock Springer New York, 2002.
\newblock ISBN 9780387215822.
\newblock \doi{10.1007/978-0-387-21582-2}.

\bibitem[Donea(1984)]{TI_Donea1984a}
J.~Donea.
\newblock A taylor--galerkin method for convective transport problems.
\newblock \emph{International Journal for Numerical Methods in Engineering},
  20\penalty0 (1):\penalty0 101--119, 1984.
\newblock \doi{10.1002/nme.1620200108}.

\bibitem[Dutt et~al.(2000)Dutt, Greengard, and Rokhlin]{TI_Dutt2000a}
A.~Dutt, L.~Greengard, and V.~Rokhlin.
\newblock Spectral deferred correction methods for ordinary differential
  equations.
\newblock \emph{Bit Numerical Mathematics}, 40\penalty0 (2):\penalty0 241--266,
  2000.
\newblock \doi{10.1023/a:1022338906936}.

\bibitem[Fehn et~al.(2017)Fehn, Wall, and Kronbichler]{TI_Fehn2017a}
N.~Fehn, W.~A. Wall, and M.~Kronbichler.
\newblock On the stability of projection methods for the incompressible
  {N}avier-{S}tokes equations based on high-order discontinuous {G}alerkin
  discretizations.
\newblock \emph{Journal of Computational Physics}, 351:\penalty0 392--421,
  2017.
\newblock \doi{10.1016/j.jcp.2017.09.031}.

\bibitem[Frank et~al.(1997)Frank, Hundsdorfer, and Verwer]{TI_Frank1997a}
J.~Frank, W.~Hundsdorfer, and J.~Verwer.
\newblock On the stability of implicit-explicit linear multistep methods.
\newblock \emph{Applied Numerical Mathematics}, 25\penalty0 (2-3):\penalty0
  193--205, 1997.
\newblock \doi{10.1016/s0168-9274(97)00059-7}.

\bibitem[Frolkovi{\v c} and {\v Z}erav{\' y}(2023)]{TI_Frolkovic2023a}
P.~Frolkovi{\v c} and M.~{\v Z}erav{\' y}.
\newblock High resolution compact implicit numerical scheme for conservation
  laws.
\newblock \emph{Applied Mathematics and Computation}, 442:\penalty0 127720,
  2023.
\newblock \doi{10.1016/j.amc.2022.127720}.

\bibitem[Glaubitz et~al.(2019)Glaubitz, Nogueira, Almeida, Cant{\~a}o, and
  Silva]{SE_Glaubitz2019a}
J.~Glaubitz, A.~C. Nogueira, J.~L.~S. Almeida, R.~F. Cant{\~a}o, and C.~A.~C.
  Silva.
\newblock Smooth and compactly supported viscous sub-cell shock capturing for
  discontinuous galerkin methods.
\newblock \emph{Journal of Scientific Computing}, 79\penalty0 (1):\penalty0
  249--272, 2019.
\newblock \doi{10.1007/s10915-018-0850-3}.

\bibitem[Golub and Ye(1999)]{KR_Golub1999a}
G.~H. Golub and Q.~Ye.
\newblock Inexact preconditioned conjugate gradient method with inner-outer
  iteration.
\newblock \emph{SIAM Journal on Scientific Computing}, 21\penalty0
  (4):\penalty0 1305--1320, 1999.

\bibitem[Gottlieb(2005)]{TI_Gottlieb2005a}
S.~Gottlieb.
\newblock On high order strong stability preserving {R}unge-{K}utta and multi
  step time discretizations.
\newblock \emph{Journal of Scientific Computing}, 25\penalty0 (1--2):\penalty0
  105--128, 2005.
\newblock \doi{10.1007/bf02728985}.

\bibitem[Gottlieb and Ketcheson(2016)]{TI_Gottlieb2016a}
S.~Gottlieb and D.~Ketcheson.
\newblock \emph{Time Discretization Techniques}, pages 549--583.
\newblock Elsevier, 2016.
\newblock \doi{10.1016/bs.hna.2016.08.001}.

\bibitem[Guesmi et~al.(2023)Guesmi, Grotteschi, and Stiller]{TI_Guesmi2023a}
M.~Guesmi, M.~Grotteschi, and J.~Stiller.
\newblock Assessment of high‐order imex methods for incompressible flow.
\newblock \emph{International Journal for Numerical Methods in Fluids},
  95\penalty0 (6):\penalty0 954--978, 2023.
\newblock \doi{10.1002/fld.5177}.

\bibitem[Hagstrom and Zhou(2006)]{TI_Hagstrom2006a}
T.~Hagstrom and R.~Zhou.
\newblock On the spectral deferred correction of splitting methods for initial
  value problems.
\newblock \emph{Communications in Applied Mathematics and Computational
  Science}, 1\penalty0 (1):\penalty0 169--205, 2006.
\newblock \doi{10.2140/camcos.2006.1.169}.

\bibitem[Hairer and Wanner(1996)]{TI_Hairer1996a}
E.~Hairer and G.~Wanner.
\newblock \emph{Solving Ordinary Differential Equations {II}}.
\newblock Springer Berlin Heidelberg, 1996.
\newblock \doi{10.1007/978-3-642-05221-7}.

\bibitem[Hairer et~al.(1993)Hairer, N{\o}rsett, and Wanner]{TI_Hairer1993a}
E.~Hairer, S.~P. N{\o}rsett, and G.~Wanner.
\newblock \emph{Solving Ordinary Differential Equations I}.
\newblock Springer Berlin Heidelberg, 1993.
\newblock \doi{10.1007/978-3-540-78862-1}.

\bibitem[Hansen and Strain(2011)]{TI_Hansen2011a}
A.~C. Hansen and J.~Strain.
\newblock On the order of deferred correction.
\newblock \emph{Applied Numerical Mathematics}, 61\penalty0 (8):\penalty0
  961--973, 2011.
\newblock \doi{10.1016/j.apnum.2011.04.001}.

\bibitem[Harten and Hyman(1983)]{FV_Harten1983a}
A.~Harten and J.~M. Hyman.
\newblock Self adjusting grid methods for one-dimensional hyperbolic
  conservation laws.
\newblock \emph{Journal of Computational Physics}, 50\penalty0 (2):\penalty0
  235--269, 1983.
\newblock \doi{10.1016/0021-9991(83)90066-9}.

\bibitem[Hestenes and Stiefel(1952)]{KR_Hestenes1952a}
M.~R. Hestenes and E.~Stiefel.
\newblock Methods of conjugate gradients for solving linear systems.
\newblock \emph{Journal of research of the National Bureau of Standards},
  49\penalty0 (6):\penalty0 409--436, 1952.

\bibitem[Hesthaven and Warburton(2008)]{SE_Hesthaven2008a}
J.~S. Hesthaven and T.~Warburton.
\newblock \emph{Nodal Discontinuous {G}alerkin Methods}.
\newblock Springer, 2008.

\bibitem[Hughes and Mallet(1986)]{FE_Hughes1986c}
T.~J.~R. Hughes and M.~Mallet.
\newblock A new finite element formulation for computational fluid dynamics:
  {III}. {T}he generalized streamline operator for multidimensional
  advective-diffusive systems.
\newblock \emph{Computer Methods in Applied Mechanics and Engineering},
  58:\penalty0 305--328, 1986.

\bibitem[Hundsdorfer and Verwer(2003)]{TI_Hundsdorfer2003a}
W.~Hundsdorfer and J.~Verwer.
\newblock \emph{Numerical Solution of Time-Dependent
  Advection-Diffusion-Reaction Equations}.
\newblock Springer Berlin Heidelberg, 2003.
\newblock \doi{10.1007/978-3-662-09017-6}.

\bibitem[Huynh(2023)]{TI_Huynh2023a}
H.~T. Huynh.
\newblock Discontinuous {G}alerkin and related methods for {ODE}.
\newblock \emph{Journal of Scientific Computing}, 96\penalty0 (2), 2023.
\newblock \doi{10.1007/s10915-023-02233-2}.

\bibitem[Izzo and Jackiewicz(2017)]{TI_Izzo2017a}
G.~Izzo and Z.~Jackiewicz.
\newblock Highly stable implicit--explicit runge--kutta methods.
\newblock \emph{Applied Numerical Mathematics}, 113:\penalty0 71--92, 2017.
\newblock \doi{10.1016/j.apnum.2016.10.018}.

\bibitem[Janssen and Kanschat(2011)]{SE_Janssen2011a}
B.~Janssen and G.~Kanschat.
\newblock Adaptive multilevel methods with local smoothing for $h^1$- and
  $h^{\mathrm{curl}}$-conforming high order finite element methods.
\newblock \emph{{SIAM} Journal on Scientific Computing}, 33\penalty0
  (4):\penalty0 2095--2114, 2011.
\newblock \doi{10.1137/090778523}.

\bibitem[Jaust et~al.(2016)Jaust, Sch{\"u}tz, and Seal]{TI_Jaust2016a}
A.~Jaust, J.~Sch{\"u}tz, and D.~C. Seal.
\newblock Implicit multistage two-derivative {D}iscontinuous {G}alerkin schemes
  for viscous conservation laws.
\newblock \emph{Journal of Scientific Computing}, 69\penalty0 (2):\penalty0
  866--891, 2016.
\newblock \doi{10.1007/s10915-016-0221-x}.

\bibitem[Karniadakis et~al.(1991)Karniadakis, Israeli, and
  Orszag]{TI_Karniadakis1991a}
G.~E. Karniadakis, M.~Israeli, and S.~A. Orszag.
\newblock High-order splitting methods for the incompressible {N}avier-{S}tokes
  equations.
\newblock \emph{Journal of Computational Physics}, 97\penalty0 (2):\penalty0
  414--443, 1991.
\newblock \doi{10.1016/0021-9991(91)90007-8}.

\bibitem[Kennedy and Carpenter(2003)]{TI_Kennedy2003a}
C.~A. Kennedy and M.~H. Carpenter.
\newblock Additive {R}unge-{K}utta schemes for convection-diffusion-reaction
  equations.
\newblock \emph{Applied Numerical Mathematics}, 44\penalty0 (1-2):\penalty0
  139--181, 2003.
\newblock \doi{10.1016/s0168-9274(02)00138-1}.

\bibitem[Kennedy and Carpenter(2016)]{TI_Kennedy2016a}
C.~A. Kennedy and M.~H. Carpenter.
\newblock Diagonally implicit runge{\textendash}kutta methods for ordinary
  differential equations: A review.
\newblock Technical Report NASA/TM 2016-219173, NASA, 2016.

\bibitem[Ketcheson(2008)]{TI_Ketcheson2008a}
D.~I. Ketcheson.
\newblock Highly efficient strong stability-preserving {R}unge--{K}utta methods
  with low-storage implementations.
\newblock \emph{SIAM Journal on Scientific Computing}, 30\penalty0
  (4):\penalty0 2113--2136, 2008.
\newblock \doi{10.1137/07070485x}.

\bibitem[Klein et~al.(2015)Klein, Kummer, Keil, and Oberlack]{TI_Klein2015a}
B.~Klein, F.~Kummer, M.~Keil, and M.~Oberlack.
\newblock An extension of the {SIMPLE} based discontinuous {G}alerkin solver to
  unsteady incompressible flows.
\newblock \emph{International Journal for Numerical Methods in Fluids},
  77\penalty0 (10):\penalty0 571--589, 2015.
\newblock \doi{10.1002/fld.3994}.

\bibitem[Kl\"{o}ckner et~al.(2011)Kl\"{o}ckner, Warburton, and
  Hesthaven]{SE_Kloeckner2011a}
A.~Kl\"{o}ckner, T.~Warburton, and J.~S. Hesthaven.
\newblock Viscous shock capturing in a time-explicit {D}iscontinuous {G}alerkin
  method.
\newblock \emph{Mathematical Modelling of Natural Phenomena}, 6\penalty0
  (3):\penalty0 57--83, 2011.
\newblock \doi{10.1051/mmnp/20116303}.

\bibitem[Lax and Wendroff(1960)]{TI_Lax1960a}
P.~Lax and B.~Wendroff.
\newblock Systems of conservation laws.
\newblock \emph{Communications on Pure and Applied Mathematics}, 13\penalty0
  (2):\penalty0 217--237, 1960.
\newblock \doi{10.1002/cpa.3160130205}.

\bibitem[Layton and Minion(2007)]{TI_Layton2007a}
A.~Layton and M.~Minion.
\newblock Implications of the choice of predictors for semi-implicit {P}icard
  integral deferred correction methods.
\newblock \emph{Communications in Applied Mathematics and Computational
  Science}, 2\penalty0 (1):\penalty0 1--34, 2007.
\newblock \doi{10.2140/camcos.2007.2.1}.

\bibitem[Layton and Minion(2005)]{TI_Layton2005a}
A.~T. Layton and M.~L. Minion.
\newblock Implications of the choice of quadrature nodes for {P}icard integral
  deferred corrections methods for ordinary differential equations.
\newblock \emph{BIT Numerical Mathematics}, 45\penalty0 (2):\penalty0 341--373,
  2005.
\newblock \doi{10.1007/s10543-005-0016-1}.

\bibitem[Lottes and Fischer(2005)]{SE_Lottes2005a}
J.~W. Lottes and P.~F. Fischer.
\newblock Hybrid multigrid/{S}chwarz algorithms for the spectral element
  method.
\newblock \emph{Journal of Scientific Computing}, 24\penalty0 (1):\penalty0
  45--78, 2005.
\newblock \doi{10.1007/s10915-004-4787-3}.

\bibitem[Macca and Boscarino(2024)]{TI_Macca2024a}
E.~Macca and S.~Boscarino.
\newblock Semi-implicit-type order-adaptive cat2 schemes for systems of balance
  laws with relaxed source term.
\newblock \emph{Communications on Applied Mathematics and Computation}, 2024.
\newblock \doi{10.1007/s42967-024-00414-w}.

\bibitem[Minion and Saye(2018)]{TI_Minion2018a}
M.~Minion and R.~Saye.
\newblock Higher-order temporal integration for the incompressible
  {Na}vier-{S}tokes equations in bounded domains.
\newblock \emph{Journal of Computational Physics}, 375:\penalty0 797--822,
  2018.
\newblock \doi{10.1016/j.jcp.2018.08.054}.

\bibitem[Minion(2003)]{TI_Minion2003b}
M.~L. Minion.
\newblock Semi-implicit spectral deferred correction methods for ordinary
  differential equations.
\newblock \emph{Communications in Mathematical Sciences}, 1\penalty0
  (3):\penalty0 471--500, 2003.
\newblock \doi{10.4310/cms.2003.v1.n3.a6}.

\bibitem[Ong and Spiteri(2020)]{TI_Ong2020a}
B.~W. Ong and R.~J. Spiteri.
\newblock Deferred correction methods for ordinary differential equations.
\newblock \emph{Journal of Scientific Computing}, 83\penalty0 (3), 2020.
\newblock \doi{10.1007/s10915-020-01235-8}.

\bibitem[Pan et~al.(2021)Pan, Yan, Peir{\'{o}}, and Sherwin]{TI_Pan2021a}
Y.~Pan, Z.-G. Yan, J.~Peir{\'{o}}, and S.~J. Sherwin.
\newblock Development of a balanced adaptive time-stepping strategy based on an
  implicit {JFNK}-{DG} compressible flow solver.
\newblock \emph{Communications on Applied Mathematics and Computation}, 2021.
\newblock \doi{10.1007/s42967-021-00138-1}.

\bibitem[Pazner and Persson(2017)]{TI_Pazner2017a}
W.~Pazner and P.-O. Persson.
\newblock Stage-parallel fully implicit runge{\textendash}kutta solvers for
  discontinuous galerkin fluid simulations.
\newblock \emph{Journal of Computational Physics}, 335:\penalty0 700--717,
  2017.
\newblock \doi{10.1016/j.jcp.2017.01.050}.

\bibitem[Persson and Peraire(2006)]{SE_Persson2006a}
P.-O. Persson and J.~Peraire.
\newblock Sub-cell shock capturing for {D}iscontinuous {G}alerkin methods.
\newblock In \emph{44th AIAA Aerospace Sciences Meeting and Exhibit}. American
  Institute of Aeronautics and Astronautics, 2006.
\newblock \doi{10.2514/6.2006-112}.

\bibitem[Safjan and Oden(1995)]{TI_Safjan1995a}
A.~Safjan and J.~Oden.
\newblock High-order {T}aylor-{G}alerkin methods for linear hyperbolic systems.
\newblock \emph{Journal of Computational Physics}, 120\penalty0 (2):\penalty0
  206--230, 1995.
\newblock \doi{10.1006/jcph.1995.1159}.

\bibitem[Sch{\"u}tz and Seal(2021)]{TI_Schuetz2021a}
J.~Sch{\"u}tz and D.~C. Seal.
\newblock An asymptotic preserving semi-implicit multiderivative solver.
\newblock \emph{Applied Numerical Mathematics}, 160:\penalty0 84--101, 2021.
\newblock \doi{10.1016/j.apnum.2020.09.004}.

\bibitem[Shewchuk(1994)]{KR_Shewchuk1994a}
J.~R. Shewchuk.
\newblock An introduction to the conjugate gradient method without the
  agonizing pain.
\newblock Technical report, Pittsburgh, PA, USA, 1994.

\bibitem[Shu and Osher(1988)]{TI_Shu1988b}
C.-W. Shu and S.~Osher.
\newblock Efficient implementation of essentially non-oscillatory
  shock-capturing schemes.
\newblock \emph{Journal of Computational Physics}, 77\penalty0 (2):\penalty0
  439--471, 1988.
\newblock \doi{10.1016/0021-9991(88)90177-5}.

\bibitem[Shu and Osher(1989)]{TI_Shu1989a}
C.-W. Shu and S.~Osher.
\newblock Efficient implementation of essentially non-oscillatory
  shock-capturing schemes, ii.
\newblock \emph{Journal of Computational Physics}, 83\penalty0 (1):\penalty0
  32--78, 1989.
\newblock \doi{10.1016/0021-9991(89)90222-2}.

\bibitem[Sod(1978)]{BM_Sod1978a}
G.~A. Sod.
\newblock A survey of several finite difference methods for systems of
  nonlinear hyperbolic conservation laws.
\newblock \emph{Journal of Computational Physics}, 27\penalty0 (1):\penalty0
  1--31, 1978.
\newblock \doi{10.1016/0021-9991(78)90023-2}.

\bibitem[Speck et~al.(2015)Speck, Ruprecht, Emmett, Minion, Bolten, and
  Krause]{TI_Speck2015a}
R.~Speck, D.~Ruprecht, M.~Emmett, M.~Minion, M.~Bolten, and R.~Krause.
\newblock A multi-level spectral deferred correction method.
\newblock \emph{BIT Numerical Mathematics}, 55\penalty0 (3):\penalty0 843--867,
  2015.
\newblock \doi{10.1007/s10543-014-0517-x}.

\bibitem[Stiller(2016{\natexlab{a}})]{SE_Stiller2016a}
J.~Stiller.
\newblock Nonuniformly weighted {S}chwarz smoothers for spectral element
  multigrid.
\newblock \emph{Journal of Scientific Computing}, 72\penalty0 (1):\penalty0
  81--96, 2016{\natexlab{a}}.
\newblock \doi{10.1007/s10915-016-0345-z}.

\bibitem[Stiller(2016{\natexlab{b}})]{SE_Stiller2016b}
J.~Stiller.
\newblock Robust multigrid for high-order discontinuous {G}alerkin methods: {A}
  fast {P}oisson solver suitable for high-aspect ratio {C}artesian grids.
\newblock \emph{Journal of Computational Physics}, 327:\penalty0 317--336,
  2016{\natexlab{b}}.
\newblock \doi{10.1016/j.jcp.2016.09.041}.

\bibitem[Stiller(2020)]{TI_Stiller2020a}
J.~Stiller.
\newblock A spectral deferred correction method for incompressible flow with
  variable viscosity.
\newblock \emph{Journal of Computational Physics}, 423:\penalty0 109840, 2020.
\newblock \doi{10.1016/j.jcp.2020.109840}.

\bibitem[Tavelli and Dumbser(2018)]{TI_Tavelli2018a}
M.~Tavelli and M.~Dumbser.
\newblock Arbitrary high order accurate space{\textendash}time discontinuous
  {G}alerkin finite element schemes on staggered unstructured meshes for linear
  elasticity.
\newblock \emph{Journal of Computational Physics}, 366:\penalty0 386--414,
  2018.
\newblock \doi{10.1016/j.jcp.2018.03.038}.

\end{thebibliography}
\end{document}